\documentclass[a4paper, 11pt]{article}

\usepackage{a4wide}
\usepackage{graphicx,array,xurl,fancyvrb,hyperref,placeins,verbatim}
\usepackage{amsmath}
\usepackage{amssymb}
\usepackage{amsthm}
\usepackage{amsfonts}
\usepackage{mathtools}

\usepackage[backend=biber,style=authoryear, url=false, doi=false, isbn=false]{biblatex}
\usepackage{tikz-cd}
\usetikzlibrary{positioning,intersections}
\usepackage{graphicx}
\usepackage[all]{xy}
\usepackage{xcolor}
\usepackage{subcaption}

\theoremstyle{plain}
\newtheorem{theorem}{Theorem}[section]
\newtheorem{lemma}[theorem]{Lemma}
\newtheorem{proposition}[theorem]{Proposition}

\theoremstyle{definition}
\newtheorem{definition}[theorem]{Definition}
\newtheorem{example}[theorem]{Example}

\theoremstyle{remark}
\newtheorem{remark}[theorem]{Remark}

\newcommand{\Z}{\mathbb{Z}}

\newcommand{\R}{\mathbb{R}}

\newcommand{\rrrarrows}{\substack{\longrightarrow\\[-0.85em] \longrightarrow \\[-0.85em] \longrightarrow}} 

\begin{document}
	
	\title{Morita equivalence of shifted symplectic Lie $n$-groupoids}
	
	\author{Milena Weiershausen \\
		Department of Mathematics\\
		University of Göttingen\\
		Bunsenstraße 3--5\\
		37085 Göttingen, Germany\\
		\texttt{a.weiershausen@stud.uni-goettingen.de}
	}
	
	\date{\today}
	
	\maketitle
	
	\begin{abstract}
		Symplectic structures on higher objects like Lie groupoids have been studied for some time now, but not all of the proposed definitions are preserved under Morita equivalence of Lie groupoids, in turn giving rise to a consistent notion of symplectic stacks. Recently, this concept has been generalized to $m$-shifted symplectic forms on Lie $n$-groupoids, which are indeed preserved under Morita equivalence of Lie $n$-groupoids. In this paper, we give a rigorous proof for this statement.
	\end{abstract}
	
	\tableofcontents

	\section{Introduction}
	Noether's theorem shows that the conservation laws of a physical system are governed by the continuous symmetries of this system (cf.~\cite{noether1918invariante}). This means that the system is invariant under the action of a specific Lie group. Given a Lie group action on a manifold, however, the orbit space in does not behave nicely in general: The quotient $M/G$ is often not a manifold itself (cf.~\cite{blohmann2008stacky}). Resolving this requires considering the differentiable stack associated to the action Lie groupoid $M\times G \rightrightarrows M$ (this is sometimes called the stacky quotient $M//G$ of this action, cf.~\cite{noohi2005foundations}). This action Lie groupoid is built by not identifying objects in the same orbit and instead viewing the equivalences between them given by the action of group elements as arrows between them. This construction remembers the Lie group action as well as the smooth structure of both the object space $M$ and the arrow space $M\times G$. Due to this kind of construction, Lie groupoids are crucial objects in gauge theories, while higher gauge theories need to make use of higher Lie groupoids (cf.~\cite{alfonsi2023higher}).  These higher objects can be expressed using the language of simplicial objects and Kan conditions.
	
	On the other hand, symplectic structures are also essential in modern physics. Even in the simplest cases, symplectic manifolds are the natural structure of the phase space in the Hamiltonian description of classical mechanics (cf.~\cite{bates1997lectures}). From there, symplectic structures arise in many different areas of physics. An interesting question now is whether there is a nice notion of symplectic structures on higher objects like Lie $n$-groupoids. When looking for symplectic structures on smooth spaces with additional algebraic structure like Lie group(oid)s and higher versions thereof, it is natural to ask the symplectic form to be compatible with the algebraic structure. This property is called multiplicativity. While there are no nontrivial multiplicative $2$-forms on Lie groups, Lie groupoids do have a notion of multiplicative symplectic forms. This definition of symplectic Lie groupoids is a special case of twisted symplectic Lie groupoids, which can be described using infinitesimal multiplicative forms on the corresponding Lie algebroids and have now been generalized to $m$-shifted symplectic Lie $n$-groupoids (cf.~\cite{cuecazhu2023shiftedsymplectic}).
	
	Many of these building blocks have only been put together at the start of this century and progress in some of the generalizations and implications is still not finished. Lie groupoids as a many-object version of Lie groups have been known since Ehresmann introduced them as differentiable groupoids (cf.~\cite{ASENS_1963_3_80_4_349_0}) and stacks have been introduced as higher sheaves by Grothendieck, 
	but even though algebraic stacks have been commonly used in algebraic geometry for a long time, it took a while before topological and differentiable stacks were introduced as their counterparts in topology and differential geometry, respectively (cf.~\cite{noohi2005foundations}). The intuition behind stacks is that instead of identifying equivalent objects, one remembers all of the equivalencies as invertible arrows acting on objects. This gives the clusters of equivalent objects the structure of a groupoid. As it turns out, stacks over smooth manifolds with additional smooth structure in the form of a differentiable atlas can be represented by Lie groupoids. Conversely, for a given Lie groupoid $G_{\bullet}$, the category of $G_{\bullet}$-principal bundles has the structure of a differentiable stack $\mathcal{B}G$. It is called the classifying stack of the Lie groupoid $G_{\bullet}$ because it has similar properties to the topological classifying space $BG$ which gives rise to the universal principal bundle $EG \rightarrow BG$ and is given in this context by the geometric realization of the simplicial nerve $NG_{\bullet}$ of $G_{\bullet}$. This construction of $\mathcal{B}G$ leads to a correspondence between Lie groupoids up to a specific kind of equivalence called Morita equivalence and differentiable stacks up to isomorphism. This correspondence is described in detail in \cite{behrend2011diffstacksgerbes}. In practice, this means that we can view Lie groupoids as models for differentiable stacks, and there can be many different (Morita equivalent) models describing the same stack. The simplicial description of Lie groupoids, which leads to the generalization to Lie $n$-groupoids, was discovered by Henriquez (cf.~\cite{henriques2008integrating}), who reframed Duskin's definition of simplicial $n$-groupoids (cf.~\cite{duskin1979higher}) from the world of sets to that of smooth manifolds. Grothendieck proved that Lie $1$-groupoids in this description are the same thing as simplicial nerves of traditional Lie groupoids (cf.~\cite{grothendieck1961techniques}).
	
	In the meantime, definitions of symplectic structures on Lie groupoids (cf.~\cite{weinstein1987symplectic}) started appearing at the same time as Lie groupoids and their corresponding Lie algebroids were gaining attention (cf.~\cite{kosmann2016multiplicativity} for an overview of these developments). A few decades later, the discussion of infinitesimal multiplicative structures and IM-forms led to a more general notion of twisted symplectic groupoids (also called quasi-symplectic groupoids) and their presymplectic versions (cf.~\cite{BCWZ2004dirac}, \cite{cattaneo2004integration}, \cite{xu2004momentum}, and \cite{bursztyn2009linear}). This was only recently generalized to Lie $n$-groupoids: Based on the very general description of shifted symplectic structures on derived $n$-stacks in the setting of algebraic geometry in \cite{PTVV2013shifted} and important contributions to the description of this theory in the world of smooth manifolds in \cite{getzlerslides}, the theory of $m$-shifted symplectic Lie $n$-groupoids was fully developed with a suitable notion of symplectic Morita equivalence in \cite{cuecazhu2023shiftedsymplectic}.
	
	The main goal of this paper is to write down a proof for the fact that shifted symplectic structures are inherited under Morita equivalence of Lie $n$-groupoids for any $n\geq 0$. This statement has been known for a while, but no proof appears in the literature yet, so we will provide a rigorous proof in Theorem (\ref{main}), which is the main theorem of the paper: 
	\begin{theorem}
		Let
		\begin{center}
			\begin{tikzcd}
				& Z_{\bullet} \arrow[swap, twoheadrightarrow]{dl}{g_{\bullet}} \arrow[twoheadrightarrow]{dr}{h_{\bullet}} & \\
				X_{\bullet} & & Y_{\bullet}
			\end{tikzcd}\\
		\end{center}
		be a Morita equivalence of Lie $n$-groupoids and $\alpha_{\bullet}$ an $m$-shifted symplectic form on $X_{\bullet}$. Then, there is also an induced $m$-shifted symplectic form $\beta_{\bullet}$ on $Y_{\bullet}$ and an $(m-1)$-shifted $2$-form $\phi_{\bullet}$ on $Z_{\bullet}$ such that the zig-zag of hypercovers becomes a symplectic Morita equivalence.
	\end{theorem}
	Leading up to this proof, we will briefly explain the simplicial notation for Lie $n$-groupoids and then describe how $m$-shifted symplectic structures have been defined on them in \cite{cuecazhu2023shiftedsymplectic}. Then, we will discuss symplectic Morita equivalence between these structures and give a detailed proof for cohomological descent of hypercovers since this will be the central property used in the proof of the main theorem. The theorem is especially useful when dealing with multiple Lie $n$-groupoid models for the same differentiable $n$-stack, because it allows us to move shifted symplectic structures from one model to another via the Morita equivalence.

	
	\section{Lie $n$-groupoids}
	
	A \emph{simplicial object} in a category $\mathcal{C}$ can be thought of as a tower of objects in that category with special structure maps between them (the face and degeneracy maps). We will be discussing simplicial objects in the categories of sets and smooth manifolds respectively, but it should be kept in mind that the constructions presented in the following can also be used for Banach manifolds in the infinite-dimensional case.
	
	A \emph{simplicial set} [resp.~\emph{simplicial manifold}] is a contravariant functor $X_{\bullet}$ from the category $\Delta$ of ordered sets to $\textbf{Set}$ [resp.~$\textbf{Mfd}$]. We denote $X_{\bullet}([m])$ as $X_m$. The behaviour of $X_{\bullet}$ on morphisms is determined by the face and degeneracy maps $d^{m}_{i} = X_{\bullet}({\delta}^{i}_{m}): X_m \rightarrow X_{m-1}$ with ${\delta}^{i}_{m}: [m-1] \rightarrow [m]$, sending $k\in \{0,\dots, i-1\}$ to itself and $k\in \{i,\dots, m\}$ to $k+1$, and ${\sigma}^{m}_{i} = X_{\bullet}(s^{i}_{m}): X_m \rightarrow X_{m+1}$ with ${\sigma}^{i}_{m}: [m+1] \rightarrow [m]$, sending $k\in \{0,\dots, i\}$ to itself and $k\in \{i+1,\dots, m\}$ to $k-1$, satisfying the following simplicial identities:
	\begin{align}
		d^{n-1}_{i} \circ d^{n}_{j} &= d^{n-1}_{j-1} \circ d^{n}_{i}\text{, if } i<j,\\
		s^{n}_{i} \circ s^{n-1}_{j} &= s^{n}_{j+1} \circ s^{n-1}_{i}\text{, if } i\leq j,\\
		d^{n}_{i} \circ s^{n-1}_{j} &= s^{n-2}_{j-1} \circ d^{n-1}_{i}\text{, if } i<j,\\
		d^{n}_{i} \circ s^{n-1}_{j} &= s^{n-2}_{j} \circ d^{n-1}_{i-1}\text{, if } i>j+1,\\
		d^{n}_{j} \circ s^{n-1}_{j} &= \mbox{id} = d^{n}_{j+1} \circ s^{n-1}_{j}.  
	\end{align}
	
	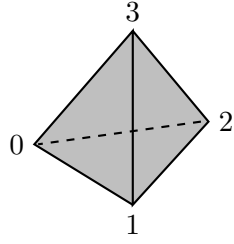
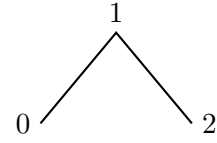
\begin{figure}[htbp]
		\centering
		
		\begin{subfigure}[b]{0.45\textwidth}
			\centering
			
			\begin{tikzpicture}
				
				\draw[thick, fill=lightgray] (0,1.6) -- (1.0,0.4) -- (0,-0.7) -- (0,1.6) -- (-1.3,0.1) -- (0,-0.7);
				\draw[thick, dashed] (0.9,0.4) -- (-1.25,0.1) ;
				
				\fill (-1.3, 0.1) circle (0pt) node[left] {0};
				\fill (0,-0.7) circle (0pt) node[below] {1};
				\fill (1,0.4) circle (0pt) node[right] {2};
				\fill (0,1.6) circle (0pt) node[above] {3};
				
			\end{tikzpicture}
			\caption{The standard simplex $\Delta^3$.}
			\label{fig1a}
		\end{subfigure}
		\begin{subfigure}[b]{0.45\textwidth}
			\centering
			
			\begin{tikzpicture}
				
				\draw[thick] (-1,0.4) -- (0,1.6) -- (1.0,0.4);
				
				\fill (-1,0.4) circle (0pt) node[left] {0};
				\fill (1,0.4) circle (0pt) node[right] {2};
				\fill (0,1.6) circle (0pt) node[above] {1};
				
			\end{tikzpicture}
			\caption{The horn $\Lambda^2_1$.}
			\label{fig1b}
			
		\end{subfigure}
		
		\caption{Standard examples of simplicial sets.}
		\label{fig1}
	\end{figure}
	
	Two standard examples we will be referring to frequently are the standard simplex and the horn. The \emph{standard $m$-simplex} $\Delta^m$ is given by $(\Delta^m)_k \coloneqq \{ f: [k] \rightarrow [m] \text{ } \vert \text{ } f(i) \leq f(j) \text{ for all }  0\leq i\leq j\leq k \}$ (see Fig.\ref{fig1a}). The \emph{horn} $\Lambda^m_j$ on the other hand is given by $(\Lambda^m_j)_k \coloneqq \{ f \in (\Delta^m)_k \vert \{0, \dots, j-1, j+1, \dots, m \} \nsubseteq \{f(0), \dots, f(k) \} \}$ (see Fig.\ref{fig1b}). In general, we can think of the horn as removing the highest non-degenerate face as well as the face opposite the $j$-th vertex from a simplex. Note that from the definition, $\Delta^m$ and $\Lambda^m_j$ are only simplicial sets, but endowed with the discrete topology, these can be viewed as simplicial manifolds. 
	
	The (higher) group(oid) structure (based on Duskin's definition of higher group(oid)s, cf. \cite{duskin1979higher}) on these simplicial manifolds is then given in the form of Kan conditions. These describe whether any horn in $X_{\bullet}$ can be (uniquely) filled to a simplex. In the diagram below, the existence of such a (unique) filling can be described in terms of the \emph{horn projection}
	\begin{equation}\label{hornproj}
		p^m_j: \mbox{Hom}(\Delta^m, X_{\bullet}) \rightarrow \mbox{Hom}(\Lambda^m_j, X_{\bullet}).
	\end{equation}
	We say that $X_{\bullet}$ satisfies the Kan condition $\mbox{Kan}(m,j)$ iff any horn $\Lambda^m_j$ in $X_{\bullet}$ can be filled to a simplex $\Delta^m$ in $X_{\bullet}$. Equivalently, $p^m_j$ is a surjective submersion. Further, we say that $X_{\bullet}$ satisfies the strict Kan condition $\mbox{Kan!}(m,j)$, iff this horn filling is unique. Equivalently, $p^m_j$ is a diffeomorphism.
	
	\begin{center}
		\begin{tikzcd}
			\Lambda^m_j \arrow{r}{} \arrow[hookrightarrow, swap]{d}{} 
			& X_{\bullet}\\
			\Delta^m \arrow[dashed, swap]{ur}{\exists lift? (*)}&
		\end{tikzcd}\\
	\end{center}
	
	\begin{definition}\label{n-gpd}
		A \emph{Lie $n$-groupoid} $X_{\bullet}$ is a simplicial manifold that satisfies:
		\begin{enumerate}
			\item $\mbox{Kan}(m,j)$ for all $m\geq 1$ and $0\leq j\leq m$.
			\item $\mbox{Kan!}(m,j)$ for all $m\geq n+1$ and $0\leq j\leq m$.
		\end{enumerate}
		A Lie $n$-group is a Lie $n$-groupoid $X_{\bullet}$ where $X_0 = \mbox{pt}$.
	\end{definition}
	
	\begin{example}\label{Ep_Liegpd}
		Given a Lie groupoid $G_1 \rightrightarrows G_0$, the simplicial nerve of this groupoid, consisting at the $k$-th level of chains of $k$ composable morphisms, has the following form:
		\begin{equation}
			\dots G_1 \times_{G_0} G_1 \rrrarrows G_1 \rightrightarrows G_0.
		\end{equation}
		This is a (simplicial) Lie $1$-groupoid together with $d_{0}(g) = s(g)$, $d_{1}(g) = t(g)$ $\forall g\in G_1$ and higher face maps consist of composing the morphisms at position $i$ and $i+1$, while the degeneracy maps simply insert a unit morphism at the $i$-th position. The Kan conditions from Definition (\ref{n-gpd}) then come from the Lie groupoid axioms of $G_{\bullet}$. If we start with a Lie group $G$ instead, we can do the same construction for the groupoid $G \rightrightarrows \mbox{pt}$. Then, we get the Lie $1$-groupoid
		\begin{equation}
			\dots G \times G \rrrarrows G \rightrightarrows \mbox{pt}.
		\end{equation}
	\end{example}
	
	We will now briefly review the suitable notion of Lie $n$-groupoid morphisms. The morphisms we are interested in will be \emph{hypercovers}. But first, note that a morphism between simplicial objects consists of levelwise morphisms commuting with face and degeneracy maps. Equivalently, one can consider the simplicial objects as functors and morphisms between them as natural transformations. Kan conditions for such a morphism $f: X_{\bullet}\rightarrow Y_{\bullet}$ mean that there exists a (unique) horn-filling for a given level in $X_{\bullet}$ in such a way that it is compatible with $f$. This is depicted in the following diagram: 
	
	\begin{center}
		\begin{tikzcd}[column sep=large, row sep=huge]
			\Lambda^m_j \arrow{r}{\varphi} \arrow[hookrightarrow]{d}{i} 
			& X_{\bullet} \arrow{d}{f} \\
			\Delta^m \arrow[dashed, swap]{ur}{\exists \phi ? (*)} \arrow{r}{\psi}
			& Y_{\bullet}
		\end{tikzcd}\\
	\end{center}
	
	We say that $f$ satisfies the Kan condition $\mbox{Kan}(m,j)$, if for any $\varphi, \psi$ as above there exists a lift $\phi$ such that $\phi \circ i = \varphi$ and $f\circ \phi = \psi$. If the lift is also unique, we say that $f$ satisfies the strict Kan condition $\mbox{Kan!}(m,j)$. Note that the unique simplicial morphism $f: X_{\bullet} \rightarrow \mbox{pt}_{\bullet}$ satisfies $\mbox{Kan}(m, j)$ iff $X_{\bullet}$ satisfies $\mbox{Kan}(m,j)$.
	
	\begin{definition}\label{Kanfib}
		A simplicial morphism $f: X_{\bullet} \rightarrow Y_{\bullet}$ between simplicial manifolds $X_{\bullet}$ and $Y_{\bullet}$ is a \emph{Kan fibration} if $f$ satisfies $\mbox{Kan}(m\geq 1, j)$ for all $0\leq j \leq m$.\\
		A simplicial morphism $f: X_{\bullet} \rightarrow Y_{\bullet}$ between Lie $n$-groupoids is a \emph{Lie $n$-groupoid Kan fibration} if $f$ satisfies $\mbox{Kan}(m\geq 1, j)$ for all $0\leq j \leq m$ and $\mbox{Kan!}(m\geq n, j)$ for all $0\leq j\leq m$.
	\end{definition}
	
	Another way to write this Kan condition for morphisms is by viewing $\mbox{Hom}(\Delta[m], X_{\bullet})$ as $X_m$ and considering
	\begin{equation}
		p_j^m: X_m \rightarrow \mbox{Hom}(\Lambda^m_j, X_{\bullet}) \times_{\mbox{Hom}(\Lambda^m_j, Y_{\bullet})} \mbox{Hom}(\Delta^m, Y_{\bullet}).
	\end{equation}
	\begin{enumerate}
		\item $\mbox{Kan}(m,j)$ holds for $f$ iff $p_j^m$ is a surjective submersion and
		\item $\mbox{Kan!}(m,j)$ holds for $f$ iff $p_j^m$ is a diffeomorphism. 
	\end{enumerate}
	Note that for $Y_{\bullet} = \mbox{pt}$, this map becomes $p^m_j : X_m \rightarrow \mbox{Hom}(\Lambda^m_j, X_{\bullet})$ from Eq.~(\ref{hornproj}) in Definition (\ref{n-gpd}), so $X_{\bullet}$ is a Lie $n$-groupoid iff $f: X_{\bullet} \rightarrow \mbox{pt}$ is a Lie $n$-groupoid Kan fibration.
	
	Definition (\ref{Kanfib}) gives us a suitable notion of morphisms between Lie $n$-groupoids. However, this concept is not quite enough to give us a notion of Morita equivalence. Hence, we need to define a stronger version of morphisms (which will correspond to \emph{acyclic} fibrations).
	
	\begin{definition}\label{hypercover}
		A simplicial morphism $f: X_{\bullet} \rightarrow Y_{\bullet}$ of Lie $n$-groupoids is a \emph{hypercover} if the maps
		\begin{equation}\label{q_m}
			q_m \coloneqq ((d_0, \dots, d_m), f_m): X_m \rightarrow \mbox{Hom}(\partial \Delta^m, X_{\bullet}) \times_{\mbox{Hom}(\partial \Delta^m, Y_{\bullet})} Y_m
		\end{equation}
		are surjective submersions for $0\leq m \leq n-1$ and a diffeomorphism for $m=n$. This automatically gives that $q_m$ is a diffeomorphism for $m > n$.
	\end{definition}
	
	Here, $\partial \Delta^m$ denotes the boundary of the simplex $\Delta^m$. The existence of the pullbacks $ \mbox{Hom}(\partial \Delta^m, X_{\bullet}) \times_{\mbox{Hom}(\partial \Delta^m, Y_{\bullet})} Y_m$ was proven in Lemma 2.4 of \cite{zhu2009stacky} using induction. Note that $q_0 = f_0$ above is a surjective submersion. In the special case $n=1$, Defintion (\ref{hypercover}) recovers the definition of a hypercover between Lie groupoids $G_{\bullet}$ and $H_{\bullet}$, as
	\begin{equation}
		q_1 = ((d_0, d_1), f_1) = ((s, t), f_1) : G_1 \rightarrow (G_0 \times G_0) \times_{H_0 \times H_0} H_1
	\end{equation}
	is an isomorphism. This means that our definition of Morita equivalence between Lie $n$-groupoids will be a generalization of Morita equivalence between Lie groupoids. When dealing with hypercovers in the later sections of this article, we will need the following proposition.
	
	\begin{proposition}\label{hypercoverfib}
		Let $f: X_{\bullet} \rightarrow Y_{\bullet}$ be a hypercover between Lie $n$-groupoids. Then, the levelwise maps $f_m: X_m \twoheadrightarrow Y_m$ are surjective submersions.
		\begin{proof}
			For $f_0: X_0 \twoheadrightarrow Y_0$, this follows immediately from the definition as $q_0 = f_0$. The general statement can then be proven using a trick from Proposition (6.5) in \cite{del2024cohomology}, which also appeared in similar forms in \cite{henriques2008integrating} and \cite{wolfson2016descent} and uses a filtration
			\begin{equation}
				S_0 = \emptyset \hookrightarrow \mbox{pt} \hookrightarrow \dots \hookrightarrow \coprod_{i=1}^{n} \mbox{pt} \hookrightarrow \dots \hookrightarrow \partial \Delta^n \hookrightarrow S_k = \Delta^n,
			\end{equation}
			where each step is the adjunction of a cell, i.e.~$S_l \cong S_{l-1} \cup_{\partial \Delta^{n_l}} \Delta^{n_l}$, to construct the pullback diagrams
			\begin{center}
				\begin{tikzcd}
					\mbox{Hom}(S_l, X_{\bullet}) \times_{\mbox{Hom}(S_l, Y_{\bullet})} Y_n \arrow[twoheadrightarrow]{r}{g_l} \arrow{d}{}
					& \mbox{Hom}(S_{l-1}, X_{\bullet}) \times_{\mbox{Hom}(S_{l-1}, Y_{\bullet})} Y_n \arrow{d}{}\\
					X_{n_l} \arrow[twoheadrightarrow]{r}{q_{n_l}} 
					& \mbox{Hom}(\partial \Delta^{n_l}, X_{\bullet}) \times_{\mbox{Hom}(\partial \Delta^{n_l}, Y_{\bullet})} Y_{n_l}
				\end{tikzcd}
			\end{center}
			Since $q_{n_l}$ from (\ref{q_m}) is a surjective submersion, so is $g_l$ and thus the composition $f_n = g_k \circ \dots \circ g_1$ is a surjective submersion.
		\end{proof}
	\end{proposition}
	
	\begin{definition}\label{ME}
		Two Lie $n$-groupoids $X_{\bullet}$ and $Y_{\bullet}$ are Morita equivalent if there exists a Lie $n$-groupoid $Z_{\bullet}$ and a zig-zag of hypercovers
		\begin{center}
			\begin{tikzcd}
				& Z_{\bullet} \arrow[swap, twoheadrightarrow]{dl}{f} \arrow[twoheadrightarrow]{dr}{g} & \\
				X_{\bullet} & & Y_{\bullet}
			\end{tikzcd}\\
		\end{center}
	\end{definition}
	
	Note that this relation is indeed an equivalence relation, as follows in particular from Proposition (\ref{equivrel}) later on.
	
	
	\section{Shifted symplectic structures}
	
	Let us now turn to the description of symplectic structures on Lie $n$-groupoids. The first components we need to define for this are tangent spaces and differential forms. But, whenever we are dealing with higher objects, we will have a tangent complex rather than a tangent space. In this case, symplectic structures can be shifted with respect to the grading of the tangent complex.
	To define the tangent complex, we will use Definition 2.8 in \cite{cuecazhu2023shiftedsymplectic}
	
	\begin{definition}\label{tangentcx}
		Let $X_{\bullet}$ be a Lie $n$-groupoid. We define the \emph{tangent complex} $(\mathcal{T}_{\bullet}X, \partial)$ as the following complex of vector bundles over $X_0$:
		\begin{equation}
			\mathcal{T}_l X \coloneqq
			\begin{cases}
				\mbox{ker} (T p^l_l) \vert_{X_0} & \text{for } l>0,\\
				TX_0 & \text{for } l=0,\\
				0 & \text{for } l<0,
			\end{cases}
		\end{equation}
		with $\partial \coloneqq (-1)^{l} Td_{l}^{l}$. Here, $p^l_l$ is the horn projection from Eq.~(\ref{hornproj}) used in Definition (\ref{n-gpd}). We write $H_{\bullet}(\mathcal{T}X)$ for the \emph{homology groups} of the tangent complex $(\mathcal{T}_{\bullet}X, \partial)$.
	\end{definition}
	
	As proven in Propositon 2.10 in \cite{cuecazhu2023shiftedsymplectic}, for $X_{\bullet}$ a Lie $n$-groupoid and $l\geq 0$, we get an isomorphism
	\begin{equation}
		\mathcal{T}_{l}X \cong TX_{l}\vert_{X_0} / \oplus_{i=0}^{l-1} \hbox{im}(Ts_{i}^{l-1}),
	\end{equation}
	sending $\partial$ to $\sum_{i=0}^{l} (-1)^{i} Td_{i}^{l}$. Additionally, we get $\mathcal{T}_{l}X = 0$ for $l>n$ since $p_{l,l}$ is an isomorphism in those cases. This tangent complex has been proven to be invariant up to weak equivalence under Morita equivalence of Lie $n$-groupoids (cf.~Corollary 2.28 in \cite{cuecazhu2023shiftedsymplectic}) and to have the structure of the corresponding Lie $n$-algebroid (cf.~\cite{li2023differentiating}), hence this is the correct notion of tangent complex for our setting. We can see this in an example for the case $n=1$.
	
	\begin{example}\label{Ep_Liealgbd}
		Let $G \rightrightarrows M$ be a Lie groupoid as in Example (\ref{Ep_Liegpd}). Then, its simplicial nerve $NG_{\bullet}$ is the Lie $1$-groupoid 
		\begin{equation}
			\dots G \times_{M} G \rrrarrows G \rightrightarrows M
		\end{equation}
		with face maps $d^2_0 = \mbox{pr}_1$, $d^2_1 = m$, and $d^2_2 = \mbox{pr}_2: G \times_{M} G \rightarrow G$ on the second level as well as $d^1_0 = s$ and $d^1_1 = t: G \rightarrow M$ on the first level. The degeneracy $s^0_0 : M \rightarrow G$ is given as $m \mapsto \mbox{id}_m$. The corresponding tangent complex from Definition (\ref{tangentcx}) thus has the form of its Lie algebroid $A = \mbox{ker}Ts\vert_{M}$,
		\begin{equation}
			\mathcal{T}_i G \coloneqq
			\begin{cases}
				A & \text{for } i=1,\\
				TM & \text{for } i=0,\\
				0 & \text{else},
			\end{cases}
		\end{equation}
		with $\partial = \rho: A \rightarrow TM$ given by the anchor $\rho = Tt\vert_{A}$.
	\end{example}
	
	An important property we will need in the proof of the main theorem of this article is the following statement (see Lemma 2.27 in \cite{cuecazhu2023shiftedsymplectic}) that hypercovers between Lie $n$-groupoids induce isomorphisms on the homology groups of the corresponding tangent complexes. 
	
	\begin{lemma}\label{hypercoverquasi}
		Let $f_{\bullet} : X_{\bullet} \rightarrow Y_{\bullet}$ be a hypercover of Lie $n$-groupoids. Then, the induced maps $Tf_l: \mathcal{T}_lX \rightarrow \mathcal{T}_lY$ form a quasi-isomorphism (i.e.~an isomorphism on the homology groups).
		\begin{proof}
			A detailed proof of this statement is provided in \cite{cuecazhu2023shiftedsymplectic}.
		\end{proof}
	\end{lemma}
	
	As done in \cite{cuecazhu2023shiftedsymplectic} and \cite{getzlerslides}, we can define the \emph{de Rham simplicial double complex} (also known as the \emph{Bott-Shulmann-Stasheff double complex}) $(\Omega^{\bullet}(X_{\bullet}), d, \delta)$ for any simplicial manifold $X_{\bullet}$, where $\Omega^{q}(X_{p})$ denotes the differential $q$-forms on the manifold $X_p$, equipped with the usual \emph{de Rham differential} $d: \Omega^{q}(X_{p}) \rightarrow \Omega^{q+1}(X_{p})$ and the \emph{simplicial differential} $\delta: \Omega^{q}(X_{p-1}) \rightarrow \Omega^{q}(X_{p})$ given as $\delta = \sum_{i=0}^{p} (-1)^{i} d_{i}^{*}$. The total complex given by $\Omega^n \coloneqq \bigoplus_{p+q=n} \Omega^{q}(X_{p})$ with differential $D = \delta + (-1)^{p} d$ gives us a notion of differential forms on simplicial manifolds. 
	
	\begin{center}
		\begin{tikzcd}
			\vdots
			& \vdots
			& \vdots & \\
			\Omega^{2}(X_{0}) \arrow{r}{\delta} \arrow{u}{d}
			& \Omega^{2}(X_{1}) \arrow{r}{\delta} \arrow{u}{d}
			& \Omega^{2}(X_{2}) \arrow{r}{\delta} \arrow{u}{d}
			& \dots \\
			\Omega^{1}(X_{0}) \arrow{r}{\delta} \arrow{u}{d}
			& \Omega^{1}(X_{1}) \arrow{r}{\delta} \arrow{u}{d}
			& \Omega^{1}(X_{2}) \arrow{r}{\delta} \arrow{u}{d}
			& \dots \\
			\Omega^{0}(X_{0}) \arrow{r}{\delta} \arrow{u}{d}
			& \Omega^{0}(X_{1}) \arrow{r}{\delta} \arrow{u}{d}
			& \Omega^{0}(X_{2}) \arrow{r}{\delta} \arrow{u}{d}
			& \dots
		\end{tikzcd}
	\end{center}
	
	To make calculations easier, we restrict to the sub-complex $\widehat{\Omega}^{\bullet}(X_{\bullet})$ of \emph{normalised differential forms}. A differential form $\alpha$ is called \emph{normalised} if it vanishes on degeneracies, i.e. 
	\begin{equation}
		\widehat{\Omega}^{\bullet}(X_{\bullet}) = \{  \alpha \in \Omega^{\bullet}(X_{\bullet}) \vert s_{i}^{*} \alpha = 0 \text{ for all (possible) } i \}.
	\end{equation}
	
	We can restrict to this subcomplex because the resulting total complexes have the same cohomology. Thus, from now on we will always consider normalised forms. Additionally, when defining (normalised) $m$-shifted $k$-forms on $X_{\bullet}$ as elements in the total complex $\widehat{\Omega}^{k+m}$, we will apply a filtration to cut off the components $\Omega^q(X_p)$ with $q<k$. As explained in \cite{getzlerslides}, this filtration was introduced in \cite{PTVV2013shifted} to be able to define closed forms in the proper way. In our case, this means that a (normalised) $m$-shifted $k$-form is a differential form $\alpha_{\bullet} \in \widehat{\Omega}^{k+m}$, which has components $\alpha_i \in \widehat{\Omega}^{k+m-i}(X_{i})$ with the highest component $\alpha_m$ being a $k$-form on $X_m$. In particular, this means that for a closed $m$-shifted $k$-form $\alpha_{\bullet}$, the last term $\alpha_m$ is \emph{multiplicative}, i.e~$\delta \alpha_m = 0$. We will use the following definition of an $m$-shifted $k$-form from \cite{cuecazhu2023shiftedsymplectic}, which matches the constructions of (twisted) symplectic groupoids in the case $n=1$:
	
	\begin{definition}\label{presympl}
		Let $X_{\bullet}$ be a simplicial manifold. A (normalised) \emph{$m$-shifted $k$-form} $\alpha_{\bullet}$ on $X_{\bullet}$ is of the form
		\begin{equation}
			\alpha_{\bullet} = \sum_{i=0}^{m} \alpha_i \text{ with } \alpha_i \in \widehat{\Omega}^{k+m-i}(X_{i}).
		\end{equation}
		Such a form $\alpha_{\bullet}$ is called \emph{closed} if $D\alpha_{\bullet} = 0$. 
		We say that a closed (normalised) $m$-shifted $2$-form $\alpha_{\bullet}$ on $X_{\bullet}$ is an \emph{$m$-shifted presymplectic form} on $X_{\bullet}$. 
	\end{definition}
	
	We can alternatively write the filtration as 
	\begin{align*}
		F^{k} \widehat{\Omega}^q(X_p) \coloneqq \begin{cases}
			0 \text{,} & \text{if } q < k,\\
			\widehat{\Omega}^q(X_p) \text{,} & \text{if } q\geq k.
		\end{cases}
	\end{align*}
	Accordingly, we denote the truncated double complex as $F^k \widehat{\Omega}^{\bullet}(X_{\bullet})$ and the corresponding total complex as $F^k \widehat{\Omega}^n$. In this notation, an $m$-shifted $k$-form is an element $\alpha_{\bullet} \in F^{k} \widehat{\Omega}^{k+m}$. Note that the differentials $d$, $\delta$, and $D$ restrict to differentials on these truncated complexes. When dealing with the cohomology groups of a Lie $n$-groupoid, we will consider the cohomology of this $k$-truncated de Rham complex $H^{\ast}_D(X_{\bullet}, F^{k}\widehat{\Omega}^{\bullet})$ for some fixed $k\geq 0$.
	
	As is expressed in the name \emph{presymplectic} form, we now have almost all of the ingredients we need to define shifted symplectic forms on Lie $n$-groupoids in Definition (\ref{presympl}). We have a concept of $m$-shifted $k$-forms, in particular $m$-shifted $2$-forms, and we can also ascertain whether such a form is closed or not. The only property we are missing is the non-degeneracy of such a form. This, however, is more difficult to define in this generalized setting. Usually, when dealing with a presymplectic form $\omega_{\bullet} \in \Omega^2(X)$ on a manifold $X$, we would consider the associated map $\omega_{\bullet}^{\flat}: T X \rightarrow T^{*}X$ between the tangent and the cotangent bundle and check whether this is an isomorphism. In our case, the $2$-form $\omega_m \in \Omega^2(X_m)$ yields a map $\omega_m^{\flat}: TX_m \rightarrow T^{*}X_m$. This only acts on the tangent bundle. What we need, however, is a quasi-isomorphism between the tangent complex $\mathcal{T}X$ and the cotangent complex $\mathcal{T}^{*}X$ given as $\mathcal{T}_l^{*}X \coloneqq (\mathcal{T}_{-l}X)^{*}$ with differential $\partial^{*} \coloneqq \partial^t$ the transpose. Equivalently (cf.~\cite{del2020morita}), we want a pairing on the tangent complex that descends to a non-degenerate pairing on the homology groups. Taking the normalisation functor to move to the Moore complex does not work here, as the functor does not commute with the pairing. For this, we need to use the Eilenberg-Zilber isomorphism to get a pairing on the Moore complex. Because this pairing then takes place on the Lie $n$-algebroid corresponding to $X_{\bullet}$ instead of the tangent spaces themselves, the resulting pairing turns out to come from something called the \emph{infinitesimal multiplicative form} (or \emph{IM-form}) $\lambda^{\omega_{\bullet}}$ corresponding to $\omega_{\bullet}$. For the general case of a Lie $n$-groupoid, we will use an explicit formula for this IM-form that has been proven to have all of the properties we need (cf. \cite{cuecazhu2023shiftedsymplectic}). But to further motivate where this construction comes from, we will first take a look at the case $n=1$, which has been studied for a long time already and where the IM-form has been a known entity.
	
	\begin{example}\label{Ep_twistedsympl}
		Let $G \rightrightarrows M$ be a Lie groupoid and let $\omega \in \Omega^{2}(G)$ be a $2$-form on $G$ and let $\phi \in \Omega^{3}(M)$ be a closed $3$-form on $M$, i.e.~$d\phi = 0$. $\omega$ is called \emph{relatively $\phi$-closed} if $d\omega = s^{*}\phi - t^{*}\phi$ (cf.~\cite{BCWZ2004dirac}) and $\omega$ is called \emph{multiplicative} if $m^{*}\omega = \mbox{pr}_1^{*}\omega + \mbox{pr}_2^{*}\omega$ as before. Considering the simplicial nerve $NG_{\bullet}$ of the Lie groupoid $G\rightrightarrows M$ as done in Example (\ref{Ep_Liegpd}), these conditions simply translate to 
		\begin{align}
			d\phi &= 0; \\
			d\omega &= \delta \phi;\\
			\delta \omega &= 0.
		\end{align}
		Hence, we get a closed $1$-shifted $2$-form $\omega_{\bullet} = \omega + \phi$ on the Lie $1$-groupoid $NG_{\bullet}$. If $\omega_{\bullet}$ is normalised, this is a $1$-shifted presymplectic form. If $\omega$ is non-degenerate, this becomes the definition of a $\phi$-twisted symplectic groupoid from \cite{BCWZ2004dirac}, which is also known as a quasisymplectic groupoid in \cite{xu2004momentum}.\\ 
		Given a relatively $\phi$-closed multiplicative $2$-form $\omega$ as above, we can construct a vector bundle map $\sigma: A \rightarrow T^{*}M$ satisfying
		\begin{align}
			&\sigma(u)(\rho(v)) = - \sigma(v)(\rho(u));\\
			&\sigma([u,v]) =  \mathcal{L}_{\rho(u)} \sigma(v) - i_{\rho(v)} d\sigma(u) + i_{\rho(v)} i_{\rho(u)} \phi 
		\end{align}
		for all $u,v \in \Gamma(A)$ (cf.~\cite{BCWZ2004dirac}). A vector bundle map satisfying these conditions is called an \emph{IM-form} (cf.~\cite{bursztyn2009linear}). The IM-form corresponding to $\omega$ is given as $\sigma(u) = i_{u}\omega\vert_{TM}$ for all $u\in A$. If $G$ is $s$-simply connected over $M$, the correspondence between relatively $\phi$-closed multiplicative $2$-forms $\omega$ and IM-forms $\sigma$ is one-to-one (cf.~\cite{BCWZ2004dirac}).\\ 
		Another way to look at this construction is via the fact that a $2$-form on a Lie groupoid $G \rightrightarrows M$ is multiplicative iff the bundle map $\omega^{\flat}: TG \rightarrow T^{*}G$ is a Lie groupoid morphism (cf.~\cite{bursztyn2009linear}). A similar statement is also true for general multiplicative $k$-forms (cf.~\cite{bursztyn2012multiplicative}). By applying the Lie functor, multiplicative $k$-forms on Lie groupoids correspond naturally to linear infinitesimal multiplicative forms $\Lambda \in \Omega^{2}(A)$ on the total space of the corresponding Lie algebroid $A$. The condition that $\Lambda$ is infinitesimal multiplicative is equivalent to saying that $\Lambda^{\flat}: TA \rightarrow T^{*}A$ is a Lie algebroid map covering $\lambda = -\sigma^{t}: TM \rightarrow A^{*}$ (cf.~\cite{bursztyn2012multiplicative}). In the case $k=2$ and if $\omega$ is relatively $\phi$-closed, $\Lambda$ can be constructed via vector bundle maps The corresponding linear infinitesimal multiplicative form $\Lambda \in \Omega^{2}(A)$ can then be written as $\Lambda = -(\sigma^{*}\omega_{can} + \rho^{*} \tau(\phi))$ where $\omega_{can}$ is the canonical symplectic form on $T^{*}M$ and $\tau(\phi) \in \Omega^{2}(TM)$ is given pointwise as $\tau(\phi)_{X} = \mbox{pr}_{M}^{*}(i_{X} \phi)$ for $X\in TM$ and the natural projection $\mbox{pr}_M : TM \rightarrow M$ (cf.~\cite{bursztyn2009linear}). In this case, we get the following maps between the tangent and the cotangent complex of $G\rightrightarrows M$ as constructed in Example (\ref{Ep_Liealgbd}):
		\begin{center}
			\begin{tikzcd}
				A \arrow[swap]{d}{\sigma} \arrow{r}{\rho} 
				& TM \arrow{d}{\lambda}\\
				T^{*}M \arrow[swap]{r}{\rho^{*}} 
				& A^{*}
			\end{tikzcd}
		\end{center} 
		For detailed proofs, see \cite{BCWZ2004dirac}, \cite{bursztyn2009linear}, and \cite{bursztyn2012multiplicative}.
	\end{example} 
	
	For the general case of an $m$-shifted presymplectic form $\omega_{\bullet}$ on a Lie $n$-groupoid $X_{\bullet}$, we can now define the corresponding IM-form $\lambda^{\omega_{\bullet}}$ as done in \cite{cuecazhu2023shiftedsymplectic} by using an explicit formula from \cite{getzlerslides}. Afterwards, we will look at some of the properties this IM-form has and use them to define the non-degeneracy condition.
	
	\begin{definition}\label{IM_form}
		Let $X_{\bullet}$ be a Lie $n$-groupoid with an $m$-shifted presymplectic form $\omega_{\bullet}$. For $l\in \Z$, $v \in (\mathcal{T}_{l}X)_{x} \subseteq T_xX_l$, and $w \in (\mathcal{T}_{m-l}X)_{x} \subseteq T_xX_{m-l}$ in the tangent complex at $x\in X_0$, we define a pointwise pairing (formula from \cite{getzlerslides}):
		\begin{equation}
			\lambda_{x}^{\omega_{\bullet}}(v,w) \coloneqq \sum_{\sigma \in \text{Shuff}_{l, m-l}} (-1)^{\sigma} \omega_{m}(T(s_{\sigma(m-1)} \dots s_{\sigma(l)})v, T(s_{\sigma(l-1)} \dots s_{\sigma(0)})w),
		\end{equation}
		where $\text{Shuff}_{l, m-l}$ denotes the $(l, m-l)$-shuffles and $(-1)^{\sigma}$ denotes the sign of the shuffle $\sigma$.
	\end{definition}
	
	\begin{remark}\label{IM_properties}
		The idea in Definition (\ref{IM_form}) is that by viewing $(\mathcal{T}_{l}X)_{x} \subseteq T_xX_l$ and $(\mathcal{T}_{m-l}X)_{x} \subseteq T_xX_{m-l}$, we can use the degeneracies to pull $v$ and $w$ up to $TX_m$, where we can apply the $2$-form $\omega_m$. This (highest term of the) IM-form concerning $\omega_m$ is enough to define a useful notion of non-degeneracy. To see this, let us quickly remark on some of the basic properties of this IM-form, which were first mentioned in \cite{getzlerslides} and proven in the Appendix of \cite{cuecazhu2023shiftedsymplectic}. From the definition, it can be checked that $\lambda^{\omega_{\bullet}}$ is graded anti-symmetric. Additionally, because $\omega_m$ is multiplicative, $\lambda^{\omega_{\bullet}}$ is \emph{infinitesimal multiplicative}, i.e.~for $u \in (\mathcal{T}_{l+1}X)$, and $w \in (\mathcal{T}_{m-l}X)$, we have
		\begin{equation}
			\lambda^{\omega_{\bullet}}(\partial u,w) + (-1)^{l+1} \lambda^{\omega_{\bullet}}(u,\partial w) = 0.
		\end{equation}
		This means that $\lambda$ respects the differential $\partial$ on the tangent complex, making it into a chain map which hence descends to the homology groups. Thus, we get a pairing
		\begin{equation}
			\lambda^{\omega_{\bullet}}(-,-) : H_l(\mathcal{T}X) \times H_{m-l}(\mathcal{T}X) \rightarrow \R_{X_0},
		\end{equation}
		where $\R_{X_0}$ denotes the trivial bundle over $X_0$. 
	\end{remark}
	
	Using this property, we will call an $m$-shifted $2$-form $\omega_{\bullet}$ \emph{non-degenerate} if this pairing on the level of homology is pointwise non-degenerate. This gives us the final building block to define shifted symplectic forms in the same way as in Definition 2.14 in \cite{cuecazhu2023shiftedsymplectic}.
	
	\begin{definition}\label{shiftedsympl}
		A pair $(X_{\bullet}, \omega_{\bullet})$ is an \emph{$m$-shifted symplectic Lie $n$-groupoid} if $X_{\bullet}$ is a Lie $n$-groupoid and $\omega_{\bullet}$ is a closed, normalised, and non-degenerate $m$-shifted $2$-form on $X_{\bullet}$.
	\end{definition}
	
	The following Examples (\ref{Ep_shiftedsymplmfd}) and (\ref{Ep_shiftedsymplgpd}) show that Definition (\ref{shiftedsympl}) recovers important special cases.
	
	\begin{example}\label{Ep_shiftedsymplmfd}
		Let $M$ be a manifold. The simplicial nerve has the form
		\begin{equation}
			\dots M \rrrarrows M \rightrightarrows M
		\end{equation}
		with all the arrows being the identity (under identification of $m\in M$ with the arrow $\mbox{id}_m$). In this case, the tangent complex is concentrated only in degree $l=0$ with $\mathcal{T}_0 M = TM$ and $\mathcal{T}_l M = 0$ for all $l \neq 0$. Thus, we can only have a $0$-shifted symplectic structure in $M$. Let $\omega_{\bullet}$ be a $0$-shifted presymplectic form on $M$. This means that $\omega_{\bullet} = \omega \in \Omega^{2}(M)$ is closed, i.e.~$d\omega = 0$ and multiplicative, which in this case simply means $\mbox{id}^{*}\omega - \mbox{id}^{*}\omega = 0$ and thus is automatically satisfied. This makes sense since the algebraic structure on the unit groupoid is trivial. The IM-form corresponding to $\omega$ is given as
		\begin{equation}
			\lambda_x^{\omega}(v,w) = \omega(v, w),
		\end{equation}
		for $v, W \in \mathcal{T}_xM = T_xM$ and $x\in M$. The non-degeneracy condition is thus equivalent to the non-degeneracy of $\omega$, so a $0$-shifted symplectic Lie $0$-groupoid in Definition (\ref{shiftedsympl}) is the same as a classical symplectic manifold.
	\end{example}
	
	\begin{example}\label{Ep_shiftedsymplgpd}
		Let $G \rightrightarrows M$ be a Lie groupoid. As before, a $1$-shifted symplectic form in Definition (\ref{shiftedsympl}) on the nerve $NG_{\bullet}$ consists of $\omega_{\bullet} = \omega + \phi$ with $\omega \in \Omega^2(G)$ and $\phi \in \Omega^3(M)$ normalised, satisfying $\delta \omega = 0$, $d\omega = \delta \phi$, and $d\phi = 0$ such that the pairing induced by the IM-form $\lambda^{\omega_{\bullet}}$ is non-degenerate. The non-degeneracy in this case is equivalent to $\mbox{ker} \omega_x \cap A_x \cap \mbox{ker} \rho_x = 0$, $\forall x\in M$. (cf. \cite{cuecazhu2023shiftedsymplectic}). In total, this recovers the definition of a \emph{twisted presymplectic Lie groupoid} in \cite{BCWZ2004dirac} also known as a \emph{quasi-symplectic Lie groupoid} in \cite{xu2004momentum} as noted before. In this case, the IM-form $\lambda^{\omega_{\bullet}}$ has the form
		\begin{align}
			\lambda^{\omega_{\bullet}}(v,a) &= \pm \omega(Ts_0(v), a), \\
			\lambda^{\omega_{\bullet}}(a,v) &= \pm \omega(a, Ts_0(v)),
		\end{align}
		$\forall v\in (\mathcal{T}_0 G_{\bullet})_x = T_x M$, $a\in (\mathcal{T}_1 G_{\bullet})_x = A_x$. The induced pairing between the tangent complex and the cotangent complex then becomes:\\
		\begin{center}
			\begin{tikzcd}
				A_x \arrow[swap]{d}{\lambda^{\omega_{\bullet}}} \arrow{r}{\rho} 
				& T_x M \arrow{d}{\lambda^{\omega_{\bullet}}}\\
				(T_xM)^{*} \arrow[swap]{r}{\rho^{*}} 
				& (A_x)^{*}
			\end{tikzcd}
		\end{center}
		This matches the construction for the $n=1$ discussed in Example (\ref{Ep_twistedsympl}). In the case $\phi = 0$, we recover the classical notion of a \emph{symplectic groupoid} as defined in \cite{weinstein1987symplectic}. Note, however, that this notion is not invariant under Morita equivalence, which is why symplectic Morita equivalence only makes sense for shifted symplectic structures, as we will see in the next section.
	\end{example}

	
	\section{Symplectic Morita equivalence}
	
	Given two $m$-shifted symplectic Lie $n$-groupoids, we can define symplectic Morita equivalence between them as done in Definition 2.31 of \cite{cuecazhu2023shiftedsymplectic}:
	
	\begin{definition}\label{symplME}
		Two $m$-shifted symplectic Lie $n$-groupoids $(X_{\bullet}, \alpha_{\bullet})$ and $(Y_{\bullet}, \beta_{\bullet})$ are \emph{symplectic Morita equivalent} if there exists another Lie $n$-groupoid $Z_{\bullet}$ with an $(m-1)$-shifted $2$-form $\phi_{\bullet}$ and hypercovers
		\begin{center}
			\begin{tikzcd}
				& Z_{\bullet} \arrow[swap, twoheadrightarrow]{dl}{f_{\bullet}} \arrow[twoheadrightarrow]{dr}{g_{\bullet}} & \\
				X_{\bullet} & & Y_{\bullet}
			\end{tikzcd}\\
		\end{center}
		satisfying $f_{\bullet}^{*}\alpha_{\bullet} - g_{\bullet}^{*}\beta_{\bullet} = D\phi_{\bullet}$.
	\end{definition}
	
	Definition (\ref{symplME}) works very nicely because it says that $X_{\bullet}$ and $Y_{\bullet}$ are Morita equivalent as Lie $n$-groupoids in the sense of Definition (\ref{ME}) and, in addition, the pullbacks of the symplectic structures only differ by a gauge transformation.
	
	\begin{proposition}\label{equivrel}
		Symplectic Morita equivalence of $m$-shifted symplectic Lie $n$-groupoids is an equivalence relation.
		\begin{proof}
			Reflexivity and symmetry are easy to check, so only transitivity needs a bit more work. A proof for this can be found as Proposition 2.33 in \cite{cuecazhu2023shiftedsymplectic}, which uses $U_{\bullet} \coloneqq V_{\bullet} \times_{X_{\bullet}} Z_{\bullet}$ with natural maps $\widetilde{f'}: U_{\bullet} \rightarrow V_{\bullet}$ and $\widetilde{g}: U_{\bullet} \rightarrow Z_{\bullet}$, making
			\begin{center}
				\begin{tikzcd}
					& U_{\bullet} \arrow[swap, twoheadrightarrow]{dl}{f\circ \widetilde{f'}} \arrow[twoheadrightarrow]{dr}{g'\circ \widetilde{g}} & \\
					W_{\bullet} & & Y_{\bullet}
				\end{tikzcd}\\
			\end{center}
			a symplectic Morita equivalence by checking that
			\begin{equation}
				(f\circ \widetilde{f'})^{*}\alpha_{\bullet} - (g'\circ \widetilde{g})^{*}\gamma_{\bullet} = D((\widetilde{f'})^{*}\phi_{\bullet} + \widetilde{g}^{*}\psi_{\bullet}).
			\end{equation}
		\end{proof}
	\end{proposition}
	
	\begin{example}\label{Ep_symplME}
		As proven in Section 2 of \cite{cuecazhu2023shiftedsymplectic}, there are two basic types of symplectic Morita equivalences such that any symplectic Morita equivalence decomposes into three symplectic Morita equivalences of these two types:
		
		\emph{Strict morphisms}: Let $f_{\bullet} : X_{\bullet} \rightarrow Y_{\bullet}$ be a hypercover of Lie $n$-groupoids and $\alpha_{\bullet}$ an $m$-shifted symplectic form on $X_{\bullet}$. Then, 
		\begin{equation}
			(Y_{\bullet}, f_{\bullet}^{*} \alpha_{\bullet}) \xleftarrow{\mbox{id}_{\bullet}} (Y_{\bullet}, 0) \xrightarrow{f_{\bullet}} (X_{\bullet}, \alpha_{\bullet})
		\end{equation}
		is a symplectic Morita equivalence. The fact that $(Y_{\bullet}, f_{\bullet}^{*} \alpha_{\bullet})$ is again an $m$-shifted symplectic Lie $n$-groupoid is proven in Lemma 2.30 of \cite{cuecazhu2023shiftedsymplectic}.
		
		\emph{Gauge transformations}: Let $(X_{\bullet}, \alpha_{\bullet})$ be an $m$-shifted symplectic Lie $n$-groupoid and $\phi_{\bullet}$ an $(m-1)$-shifted $2$-form on $X_{\bullet}$. Then, $(X_{\bullet}, \alpha_{\bullet} + D\phi_{\bullet})$ is again an $m$-shifted symplectic Lie $n$-groupoid and
		\begin{equation}
			(X_{\bullet}, \alpha_{\bullet} + D\phi_{\bullet}) \xleftarrow{\mbox{id}_{\bullet}} (X_{\bullet}, \phi_{\bullet}) \xrightarrow{\mbox{id}_{\bullet}} (X_{\bullet}, \alpha_{\bullet})
		\end{equation}
		is a symplectic Morita equivalence.
	\end{example}
	
	Another important property that justifies this definition for symplectic Morita equivalence is that shifted symplectic structures are inherited under the usual Morita equivalence between Lie $n$-groupoids. This means that given a Morita equivalence 
	\begin{center}
		\begin{tikzcd}
			& Z_{\bullet} \arrow[swap, twoheadrightarrow]{dl}{g_{\bullet}} \arrow[twoheadrightarrow]{dr}{h_{\bullet}} & \\
			X_{\bullet} & & Y_{\bullet}
		\end{tikzcd}\\
	\end{center}
	between Lie $n$-groupoids $X_{\bullet}$ and $Y_{\bullet}$ and supposing that we have an $m$-shifted symplectic form $\alpha_{\bullet}$ on $X_{\bullet}$, we then get an inherited $m$-shifted symplectic form on $Y_{\bullet}$ and an $(m-1)$-shifted form on $Z_{\bullet}$ such that the Morita equivalence becomes a symplectic Morita equivalence. This property is known, but has not been explicitly written down in the literature yet. Since it will be very useful when trying to find new shifted symplectic Lie $n$-groupoid models for the same (higher) stack, this will be Theorem (\ref{main}), the main theorem in this paper. For this, we first need to prove that hypercovers between Lie $n$-groupoids satisfy cohomological descent. This statement is also known, but most authors only give sketches of a proof or only prove special cases, so we will provide a more detailed version of those proofs here.
	
	\subsection{Cohomological descent for hypercovers}
	
	The first step in proving cohomological descent for hypercovers consists of proving the result for any surjective submersion between manifolds. For this, let us quickly note how the trivial groupoid $Y_{\bullet} = (Y\rightrightarrows Y)$ behaves under cohomology as a simplicial object. In this groupoid, the objects are elements $y\in Y$ and the morphisms are the unit arrows $\mbox{id}_y$ $\forall y\in Y$. We will call this Lie groupoid the \emph{unit groupoid} corresponding to $Y$. This is a special case of the \v{C}ech cover discussed in Chapter 8 of \cite{bott1982differential} for $\mathcal{U}=\{Y \}$ and therefore yields that for any fixed $k\geq 0$, the $k$-truncated total cohomology of $Y_{\bullet}$ is the $k$-truncated de Rham cohomology of $Y$, i.e.~$H_D(Y_{\bullet}, F^k \widehat{\Omega}^{\bullet}) = H_{dR}(Y, F^k \widehat{\Omega}^{\bullet})$. In this case, the simplicial nerve $NY_{\bullet}$ has the form:
	\begin{equation}
		... Y \rrrarrows Y \rightrightarrows Y
	\end{equation}
	where all the face maps are identity (viewed as $\mbox{id}: \mbox{id}_y \mapsto y$).\\
	Thus, for fixed $q \in \Z_{\geq 0}$, we can form the complex in the simplicial direction as
	\begin{equation}
		0 \rightarrow F^k \widehat{\Omega}^q(Y) \xrightarrow{0} F^k \widehat{\Omega}^q(Y) \xrightarrow{\mbox{id}} F^k \widehat{\Omega}^q(Y) \xrightarrow{0} F^k \widehat{\Omega}^q(Y) \xrightarrow{\mbox{id}} \dots ,
	\end{equation}
	where the differential $\delta$ turns out to switch between the zero map and the identity. The cohomology ${H_{\delta}}^{l}(Y, F^k \widehat{\Omega}^q)$ of this complex is zero for $l > 0$ and $H^{0}(Y, F^k \widehat{\Omega}^q) = F^k \widehat{\Omega}^q(Y)$.
	
	\begin{lemma}\label{descent_nerves}
		Let $f: X\twoheadrightarrow Y$ be a surjective submersion between manifolds and define the Lie groupoids $X_{\bullet} = (X\times_{Y} X \rightrightarrows X)$ as the nerve groupoid corresponding to $f$ and $Y_{\bullet} = (Y \rightrightarrows Y)$ as the unit groupoid. Then, for any fixed $k\geq 0$, the map $f^{*}: F^k \widehat{\Omega}^q(Y) \hookrightarrow F^k \widehat{\Omega}^q(X)$ induces an isomorphism between the $k$-truncated de Rham cohomology of $Y$ and the $k$-truncated total cohomology of $X_{\bullet}$.
		\begin{proof}
			In the case where $f$ corresponds to a \v{C}ech cover (Case 1), this statement is proven in chapter 8 of \cite{bott1982differential} as the \emph{generalized Mayer-Vietoris principle}. The other special case to consider is the case where $f$ has a section, i.e.~a map $\sigma: Y\rightarrow X$ with $f\circ \sigma = \mbox{id}_Y$ (Case 2). We will prove this case by following \cite{getzlerslides} and considering the simplicial nerve of $X_{\bullet}$ to show that the augmented complex 
			\begin{center}
				\begin{tikzcd}
					
					& \vdots
					& \vdots
					& \vdots & \\
					0 \arrow{r}{}
					& F^k \widehat{\Omega}^{2}(Y) \arrow[hookrightarrow]{r}{f^{*}} \arrow{u}{d}
					& F^k \widehat{\Omega}^{2}(X) \arrow{r}{\delta} \arrow{u}{d}
					& F^k \widehat{\Omega}^{2}(X\times_Y X) \arrow{r}{\delta} \arrow{u}{d}
					& \dots \\
					0 \arrow{r}{}
					& F^k \widehat{\Omega}^{1}(Y) \arrow[hookrightarrow]{r}{f^{*}} \arrow{u}{d}
					& F^k \widehat{\Omega}^{1}(X) \arrow{r}{\delta} \arrow{u}{d}
					& F^k \widehat{\Omega}^{1}(X\times_Y X) \arrow{r}{\delta} \arrow{u}{d}
					& \dots \\
					0 \arrow{r}{}
					& F^k \widehat{\Omega}^{0}(Y) \arrow[hookrightarrow]{r}{f^{*}} \arrow{u}{d}
					& F^k \widehat{\Omega}^{0}(X) \arrow{r}{\delta} \arrow{u}{d}
					& F^k \widehat{\Omega}^{0}(X\times_Y X) \arrow{r}{\delta} \arrow{u}{d}
					& \dots
				\end{tikzcd}
			\end{center}
			has contractible rows using an extra codegeneracy. An \emph{extra codegeneracy} on an augmented cosimplicial cochain complex $V^{\bullet}$ is a map $s_{-1}: V^{n+1} \rightarrow V^n$ for $n \geq 0$ such that 
			\begin{align}
				s_{-1} \circ d_0 &= \mbox{id}: V^n \rightarrow V^n \text{ for } n\geq -1,\\
				s_{-1} \circ d_i &= d_{i-1} \circ s_{-1}: V^n \rightarrow V^n \text{ for } i>0,\\
				s_{-1} \circ s_i &= s_{i-1} \circ s_{-1}: V^{n+1} \rightarrow V^{n-1} \text{ for } 0\leq i \leq n.
			\end{align}
			This map $s_{-1}$ will then induce a homotopy operator 
			\begin{align}
				s_{-1} \circ \delta + \delta \circ s_{-1} & = \sum_i (-1)^{i} (s_{-1} \circ d_i + d_i \circ s_{-1})\\
				&= s_{-1} \circ d_0 + d_0 \circ s_{-1} + \sum_{i>0} (-1)^{i} (s_{-1} \circ d_i + d_i \circ s_{-1}) \\
				&= \mbox{id} + d_0 \circ s_{-1} + \sum_{i>0} (-1)^{i} (d_{i-1} \circ s_{-1} + d_i \circ s_{-1}) \\
				&= \mbox{id} + 0 = \mbox{id}
			\end{align} 
			on the rows of the double cochain complex in the figure above. In our case, this extra codegeneracy is induced by the section $\sigma: Y \rightarrow X$ of $f$ and given as $s_{-1} \coloneqq \sigma_{-1}^{*}$ for
			\begin{equation}
				\sigma_{-1} : N_{n}(X/Y) \cong Y \times_Y N_{n}(X/Y) \xrightarrow{\sigma \times \mbox{id}} X\times_Y N_{n}(X/Y) \cong N_{n+1}(X/Y).
			\end{equation} 
			The underlying map here sends $(x_0, \dots, x_n) \in N_{n}(X/Y)$ with 
			\begin{equation}
				f(x_0) = f(x_1) = \dots = f(x_n) = y
			\end{equation}
			to $(\sigma(y), x_0, \dots, x_n) \in N_{n+1}(X/Y)$. On this level, one can easily check the corresponding axioms for an extra \emph{degeneracy} (i.e.~the same axioms as before with switched composition):
			\begin{enumerate}
				\item $d_0 \circ \sigma_{-1}= \mbox{id}$:
				\begin{align}
					(d_0 \circ \sigma_{-1})(x_0, \dots, x_n) &= d_0(\sigma(y), x_0, \dots, x_n)\\
					&= (x_0, \dots, x_n)
				\end{align}
				\item $d_i \circ \sigma_{-1} = \sigma_{-1}\circ d_{i-1}$ for $i>0$:
				\begin{align}
					(d_i \circ \sigma_{-1})(x_0, \dots, x_n) &= d_i(\sigma(y), x_0, \dots, x_n)\\
					&= (\sigma(y), x_0, \dots, \widehat{x_{i-1}}, x_i, \dots, x_n)\\
					(\sigma_{-1}\circ d_{i-1})(x_0, \dots, x_n) &= \sigma_{-1}(x_0, \dots, \widehat{x_{i-1}}, x_i, \dots, x_n)\\
					&= (\sigma(y), x_0, \dots, \widehat{x_{i-1}}, x_i, \dots, x_n)
				\end{align}
				\item $s_i \circ \sigma_{-1} = \sigma_{-1} \circ s_{i-1}$ for $0\leq i \leq n$:
				\begin{align}
					(s_i \circ \sigma_{-1})(x_0, \dots, x_n) &= s_i(\sigma(y), x_0, \dots, x_n)\\
					&= (\sigma(y), x_0, \dots, x_{i-1}, x_{i-1}, x_i, \dots, x_n)\\
					(\sigma_{-1} \circ s_{i-1})(x_0, \dots, x_n) &= \sigma_{-1}(x_0, \dots, x_{i-1}, x_{i-1}, x_i, \dots, x_n)\\
					&= (\sigma(y), x_0, \dots, x_{i-1}, x_{i-1}, x_i, \dots, x_n).
				\end{align}
			\end{enumerate}
			This in turn shows that $s_{-1} = \sigma_{-1}^{*}$ is an extra codegeneracy and thus a homotopy operator on the horizontal cochain complex $F^k \widehat{\Omega}^{q}(N_{\bullet}(X/Y))$ for each $q\in \Z_{\geq 0}$. We can thus apply (\ref{homotopyop}) to conclude that the horizontal cohomology of $F^k \widehat{\Omega}^{\bullet}(N_{\bullet}(X/Y))$ vanishes for degree $>0$ and thus  $f^{*}: F^k \widehat{\Omega}^q(Y) \hookrightarrow F^k \widehat{\Omega}^q(X_0)$ induces an isomorphism between the $k$-truncated de Rham cohomology of $Y$ and the $k$-truncated total cohomology of $X_{\bullet}$ by (\ref{homotopyrmk}).
			
			To prove the statement in the general case, we will follow the proof for Proposition 2 in \cite{behrend2004cohomology}. Let $f: X\rightarrow Y$ be a surjective submersion between manifolds and let $\mathcal{V} = \{ V_i \}$ be an open cover of $Y$ such that $f$ has local sections over $\mathcal{V}$ and denote the corresponding \v{C}ech groupoid as $V_{\bullet}$. Then define the bisimplicial object $W_{\bullet \bullet}$ as $W_{mn} \coloneqq X_m \times_{Y} V_n$. We get the following diagram:
			\begin{center}
				\begin{tikzcd}
					& \vdots \arrow{d}{} \arrow[shift right]{d}{} \arrow[shift left]{d}{}
					& \vdots \arrow{d}{} \arrow[shift right]{d}{} \arrow[shift left]{d}{}
					& \vdots \arrow{d}{} \arrow[shift right]{d}{} \arrow[shift left]{d}{}\\
					\dots \arrow{r}{} \arrow[shift right]{r}{} \arrow[shift left]{r}{} 
					& W_{11} \arrow[shift left]{r}{} \arrow[shift right]{r}{} \arrow[shift left]{d}{} \arrow[shift right]{d}{}
					& W_{01} \arrow[swap, twoheadrightarrow]{r}{f_1} \arrow[shift left]{d}{} \arrow[shift right]{d}{}
					& V_1 \arrow[shift left]{d}{} \arrow[shift right]{d}{} \arrow[bend right = 35]{l}\\
					\dots \arrow{r}{} \arrow[shift right]{r}{} \arrow[shift left]{r}{} 
					& W_{10} \arrow[shift left]{r}{} \arrow[shift right]{r}{} \arrow[twoheadrightarrow]{d}{p_1}
					& W_{00} \arrow[swap, twoheadrightarrow]{r}{f_0} \arrow[twoheadrightarrow]{d}{p_0}
					& V_0 \arrow[twoheadrightarrow]{d}{p} \arrow[bend right = 35]{l}\\
					\dots \arrow{r}{} \arrow[shift right]{r}{} \arrow[shift left]{r}{} 
					& X_1 \arrow[shift left]{r}{} \arrow[shift right]{r}{} 
					& X_0 \arrow[swap, twoheadrightarrow]{r}{f} 
					& Y
				\end{tikzcd}
			\end{center}
			Note that the rows of this twice augmented bisimplicial complex are the nerve of the groupoid corresponding to $W_{0n} \twoheadrightarrow V_n$, which is a surjective submersion that admits a section. This means that the rows of the corresponding cochain complex obtained by applying the functor $F^k \widehat{\Omega}^{q}$ to everything are exact due to case 2 above. Further note that the columns of this complex are the nerve of the groupoid corresponding to $W_{m0} \twoheadrightarrow X_m$, which is a surjective submersion coming from an open cover, hence the columns of the corresponding cochain complex are exact due to case 1 above. Again using the argument about augmented double complexes with exact rows from \cite{bott1982differential}, this proves that we have the following chain of isomorphisms for fixed $q\in \Z_{\geq 0}$ on the level of cohomology
			\begin{center}
				\begin{tikzcd}
					H_D^{l}(W_{\bullet \bullet}, F^k \widehat{\Omega}^q) \arrow[swap]{d}{\cong} \arrow{r}{\cong} 
					& H_{\delta}^l(V_{\bullet}, F^k \widehat{\Omega}^q) \arrow[swap]{d}{\cong}\\
					H_{\delta}^l(X_{\bullet}, F^k \widehat{\Omega}^q)
					& H_{\delta}^l(Y\rightrightarrows Y, F^k \widehat{\Omega}^q)
				\end{tikzcd}
			\end{center}
			We can now use the argument for augmented triple complexes with exact rows from (\ref{triplecx}), which proves that $f^{*}$ induces an isomorphism on the $k$-truncated total cohomology groups.
		\end{proof}
	\end{lemma}
	
	Now that we have proven these prerequisites, we can finally turn to the general setting of a hypercover between Lie $n$-groupoids. To prove this general case, we will need a coskeleton construction used in \cite{getzlerslides} and \cite{dugger_hollander_isaksen_2004}. The details of the construction and the useful properties are based on \cite{duskin2002simplicial}.
	
	\begin{definition}\label{cosk}
		Define the \emph{$m$-skeleton} of the $n$-simplex as the simplicial set $\mbox{sk}_m \Delta^{n}$ consisting of all the $m$-simplices in $\Delta^{n}$ with degeneracies. In particular, $\mbox{sk}_{n-1} \Delta^{n} = \partial \Delta^{n}$ and $\mbox{sk}_{m} \Delta^{n} = \Delta^{n}$ for $m\geq n$. Further, define the \emph{$m$-coskeleton} of $X_{\bullet}$ as the simplicial set $\mbox{cosk}_m(X_{\bullet})$ given by $(\mbox{cosk}_m X_{\bullet})_n = \mbox{Hom}(\mbox{sk}_m \Delta^{n}, X_{\bullet})$.
	\end{definition}
	
	\begin{remark}\label{cosk_properties}
		For $n\leq m$, we get that
		\begin{align}
			(\mbox{cosk}_m X_{\bullet})_n &= \mbox{Hom}(\mbox{sk}_m \Delta^{n}, X_{\bullet})\\
			&= \mbox{Hom}(\Delta^{n}, X_{\bullet})\\
			&= X_n
		\end{align}
		and for $X_{\bullet}$ a Lie $n$-groupoid, we have that $\mbox{cosk}_m(X_{\bullet}) \cong X_{\bullet}$ for $m> n$ as simplicial sets.
	\end{remark}
	
	To deal with the coskeleta from Definition (\ref{cosk}), which are only simplicial sets, as Lie $n$-groupoids, we will need the \emph{relative coskeleton} of a hypercover as given in \cite{getzlerslides}.
	
	\begin{definition}\label{rel_cosk}
		Let $f: X_{\bullet} \twoheadrightarrow Y_{\bullet}$ be a hypercover of Lie $n$-groupoids. The \emph{relative coskeleton} corresponding to $f$ is defined as the simplicial set given by
		\begin{align} 
			(\mbox{cosk}_m(X_{\bullet}/Y_{\bullet}))_n &= (\mbox{cosk}_m(X_{\bullet}))_n \times_{(\mbox{cosk}_m(Y_{\bullet}))_n} Y_n\\
			&= \mbox{Hom}(\mbox{sk}_m \Delta^n, X_{\bullet}) \times_{\mbox{Hom}(\mbox{sk}_m \Delta^{n}, Y_{\bullet})} Y_n.
		\end{align}
		These relative coskeleta are in fact simplicial manifolds again, as proven in Lemma 2.4 of \cite{zhu2009stacky}, and have the following properties:
		\begin{enumerate}
			\item For $m > n$, we have $\mbox{cosk}_m(X_{\bullet}/Y_{\bullet}) \cong X_{\bullet}$.
			\item We can define $\mbox{cosk}_{-1}(X_{\bullet}/Y_{\bullet})$ after augmenting the simplex category by an initial object $[-1] \coloneqq \emptyset$, which yields $(\mbox{cosk}_{-1}(X_{\bullet}/Y_{\bullet}))_n = \mbox{Hom}(\emptyset, X_{\bullet}) \times_{\mbox{Hom}(\emptyset, Y_{\bullet})} Y_n = \mbox{pt} \times_{\mbox{pt}} Y_n$ and thus $\mbox{cosk}_{-1}(X_{\bullet}/Y_{\bullet}) \cong Y_{\bullet}$.
			\item We can factor the hypercover $f$ as the following tower of simplicial maps between coskeleta, which are again hypercovers (as proven in \cite{getzlerslides}).
			\begin{align}
				&X_{k} &\cong& &(\mbox{cosk}_{n+1}(X_{\bullet}/Y_{\bullet}))_k &\twoheadrightarrow& &\dots &\twoheadrightarrow& &(\mbox{cosk}_{-1}(X_{\bullet}/Y_{\bullet}))_k &\cong& &Y_{k} & \\
				&x &\mapsto& &(x, f_k(x)) &\mapsto& &\dots &\mapsto& &({i_0}^{*}\dots{i_n}^{*}x, f_k(x)) &\mapsto& &f_k(x),&
			\end{align}
			Here, $i_m : \mbox{sk}_{m-1} \Delta^k \hookrightarrow \mbox{sk}_m \Delta^k$ denotes the inclusion map.
		\end{enumerate}
	\end{definition}
	
	With this definition and the basic properties from Remark (\ref{cosk_properties}) and Definition (\ref{rel_cosk}), we can now continue with the proof of cohomological descent for a general hypercover between Lie $n$-groupoids. Note that due to the tower of hypercovers between coskeleta mentioned above, we only need to prove cohomological descent for an arbitrary one of these maps. This will become very useful in the last step of the proof.
	
	\begin{lemma}\label{descent_hypercovers}
		Hypercovers between Lie $n$-groupoids satisfy cohomological descent. This means that given a hypercover $f: X_{\bullet} \twoheadrightarrow Y_{\bullet}$ between Lie $n$-groupoids and a fixed $k\geq 0$, the corresponding map $f^{*}: F^k \widehat{\Omega}^{\bullet}(Y_{\bullet}) \rightarrow F^k \widehat{\Omega}^{\bullet}(X_{\bullet})$ is a quasi-isomorphism on the $k$-truncated total complexes, i.e.~an isomorphism on the level of the $k$-truncated total cohomology groups.
		\begin{proof} This proof will follow the one given in Proposition A.4 of \cite{dugger_hollander_isaksen_2004} and the adaptation of that proof to this setting as done in \cite{getzlerslides}. Let $f: X_{\bullet} \twoheadrightarrow Y_{\bullet}$ be a hypercover of Lie $n$-groupoids. Consider the bisimplicial manifold $Z_{ij} \coloneqq N_{i}(X_j/Y_j)$ with vertical face and degeneracy maps denoted by $d_p^{i}$ and $s_p^{i}$ respectively and with horizontal face and degeneracy maps denoted by $\bar{d}_q^{j}$ and $\bar{s}_q^{j}$, respectively. This is well-defined since the level-wise map $f_j: X_j \twoheadrightarrow Y_j$ is a surjective submersion for each level $j$ by Proposition (\ref{hypercoverfib}) and all of the structure maps are given by the nerve construction. At the lowest level, $Z_{0j} = X_j$ admits a surjective submersion $f_j$ to $Y_j$, so we can augment the picture by $Y_{\bullet}$:
			\begin{center}
				\begin{tikzcd}[cells={text width={width("$X_0 \times_{Y_0} X_0$")},align=center},
					column sep=2em]
					\hspace{-3cm} &\vdots \arrow{d}{} \arrow[shift right]{d}{} \arrow[shift left]{d}{}
					& \vdots \arrow{d}{} \arrow[shift right]{d}{} \arrow[shift left]{d}{}
					& \vdots \arrow{d}{} \arrow[shift right]{d}{} \arrow[shift left]{d}{}
					& \hspace{3cm}\\
					\dots \hspace{-3cm} 
					& X_2 \times_{Y_2} X_2 \arrow{r}{} \arrow[shift right]{r}{} \arrow[shift left]{r}{} \arrow[shift left]{d}{} \arrow[shift right]{d}{}
					& X_1 \times_{Y_1} X_1 \arrow[shift left]{r}{} \arrow[shift right]{r}{} \arrow[shift left]{d}{} \arrow[shift right]{d}{}
					& X_0 \times_{Y_0} X_0 \arrow[shift left]{d}{} \arrow[shift right]{d}{}
					&\\
					\dots \hspace{-3cm} 
					& X_2 \arrow{r}{} \arrow[shift right]{r}{} \arrow[shift left]{r}{} \arrow[twoheadrightarrow]{d}{f_2}
					& X_1 \arrow[shift left]{r}{} \arrow[shift right]{r}{} \arrow[twoheadrightarrow]{d}{f_1}
					& X_0 \arrow[twoheadrightarrow]{d}{f_0}
					&\\
					\dots \hspace{-3cm} 
					& Y_2 \arrow{r}{} \arrow[shift right]{r}{} \arrow[shift left]{r}{} 
					& Y_1 \arrow[shift left]{r}{} \arrow[shift right]{r}{} 
					& Y_0 &
				\end{tikzcd}
			\end{center}
			We can apply the de Rham functor $F^k \widehat{\Omega}^{\bullet}$ to this to obtain a triple complex given as 
			\begin{equation}
				V_{ij}^{l} \coloneqq F^k \widehat{\Omega}^{l}(Z_{ij}) \coloneqq \{\omega \in F^k \Omega^{l}(Z_{ij}) \; \vert \;  s_p \omega = 0 \text{ for } 0\leq p \leq i-1 \text{ and } \bar{s}_q \omega = 0 \text{ for } 0\leq q \leq j-1\}.
			\end{equation} 
			The two simplicial differentials on this triple complex will be denoted by 
			\begin{equation}
				\delta: V_{ij}^{l} \rightarrow V_{(i+1)j}^{l}
			\end{equation} 
			for the vertical direction and 
			\begin{equation}
				\bar{\delta}: V_{ij}^{l} \rightarrow V_{i(j+1)}^{l}
			\end{equation} 
			for the horizontal direction. We can then take the total complex $V^n \coloneqq \bigoplus_{i+j+l=n} V_{ij}^{l}$ of this to compute its cohomology. The total differential on this is given by 
			\begin{equation}
				D = \delta + (-1)^i \bar{\delta} + (-1)^{i+j}d.
			\end{equation}
			Note that $f^{*}$ descends to the subcomplexes of normalised forms because $f$ is a simplicial map and thus commutes with face and degeneracy maps.
			\begin{center}
				\begin{tikzcd}
					&\vdots 
					& \vdots 
					& \vdots 
					&\\
					& V_{1, 0}^{l} \arrow{r}{\bar{\delta}} \arrow{u}{\delta}
					& V_{1, 1}^{l} \arrow{r}{\bar{\delta}} \arrow{u}{\delta}
					& V_{1, 2}^{l} \arrow{r}{\bar{\delta}} \arrow{u}{\delta}
					& \dots \\
					& V_{0, 0}^{l} \arrow{r}{\bar{\delta}} \arrow{u}{\delta}
					& V_{0, 1}^{l} \arrow{r}{\bar{\delta}} \arrow{u}{\delta}
					& V_{0, 2}^{l} \arrow{r}{\bar{\delta}} \arrow{u}{\delta}
					& \dots \\
					& F^k \Omega^{l}(Y_0) \arrow{r}{\delta} \arrow{u}{f^{*}}
					& F^k \Omega^{l}(Y_1) \arrow{r}{\delta} \arrow{u}{f^{*}}
					& F^k \Omega^{l}(Y_2) \arrow{r}{\delta} \arrow{u}{f^{*}}
					& \dots
				\end{tikzcd}
			\end{center}
			We will now use the construction above to prove that $f^{*}: F^k \widehat{\Omega}^{\bullet}(Y_{\bullet}) \rightarrow F^k \widehat{\Omega}^{\bullet}(X_{\bullet})$ is a quasi-isomorphism between the total complexes of $F^k \widehat{\Omega}^{\bullet}(Y_{\bullet})$ and $F^k \widehat{\Omega}^{\bullet}(X_{\bullet})$ in three steps:
			
			\emph{Step 1:} Note that due to cohomological descent for nerve constructions in the Lemma (\ref{descent_nerves}), we have already proven that $f^{*}: F^k \widehat{\Omega}^{l}(Y_{j}) \rightarrow F^k \widehat{\Omega}^{l}(X_{j}) = V_{0j}^{l}$ induces a quasi-isomorphism between the $k$-truncated total complexes of differential forms on $Y_{\bullet}$ and $Z_{\bullet \bullet}$.
			
			\emph{Step 2:} Denote the \emph{diagonal complex} of $Z_{\bullet \bullet}$ as $\widehat{Z}_{\bullet}$. The simplicial differential on the double complex $\widehat{V}_{n}^{l} \coloneqq F^k \widehat{\Omega}^l (Z_{nn})$ corresponding to this diagonal is given by 
			\begin{equation}
				\widehat{\delta} = \delta \circ \bar{\delta}: \widehat{V}_{n}^{l} \rightarrow \widehat{V}_{n+1}^{l}
			\end{equation} 
			and the total differential will be denoted as $\widehat{D}$. By the Eilenberg-Zilber theorem (cf. Theorem 2.1 in \cite{eilenberg1954groups} for the original statement and Theorem 2.9 in \cite{dold1961homologie} for this generalized version), the map sending 
			\begin{equation}
				\omega \in V_{i j}^{l} \mapsto (d_{i+1}^{*})^j (\bar{d}_{0}^{*})^i \omega \in \widehat{V}_{i+j}^{l}
			\end{equation} 
			induces a quasi-isomorphism between the total complexes of $V_{\bullet \bullet}^{\bullet}$ and $\widehat{V}_{\bullet}^{\bullet}$.
			
			\emph{Step 3:} The last remaining step is to use these quasi-isomorphisms to conclude the proof that $f^{*}: F^k \widehat{\Omega}^{\bullet}(Y_{\bullet}) \rightarrow F^k \widehat{\Omega}^{\bullet}(X_{\bullet})$ is a quasi-isomorphism between the total complexes of $F^k \widehat{\Omega}^{\bullet}(Y_{\bullet})$ and $F^k \widehat{\Omega}^{\bullet}(X_{\bullet})$. Here, we will use the factorization of $f$ as hypercovers between relative coskeleta from Definition (\ref{rel_cosk}) and only show the statement for $\widetilde{f}: \mbox{cosk}_m(X_{\bullet}/Y_{\bullet}) \rightarrow \mbox{cosk}_{m-1}(X_{\bullet}/Y_{\bullet})$. As explained before, we can always factor a hypercover $f$ between Lie $n$-groupoids as a composition of these coskeleton maps, so it is enough to prove the quasi-isomorphism for an arbitrary one of these maps. Hence, for a hypercover $f: X_{\bullet} \twoheadrightarrow Y_{\bullet}$, consider an arbitrary one of the coskeleton maps 
			\begin{equation}
				\widetilde{f}: U_{\bullet} \coloneqq \mbox{cosk}_m(X_{\bullet}/Y_{\bullet}) \rightarrow \mbox{cosk}_{m-1}(X_{\bullet}/Y_{\bullet}) \eqqcolon V_{\bullet}. 
			\end{equation}
			Note that $\widetilde{f}$ is again a hypercover of Lie $n$-groupoids. Now consider the corresponding bisimplicial object $W_{i j} = N_{i}(U_{j}/V_{j})$ and its diagonal $\widehat{W}_{l}$ as before. By applying the Eilenberg-Zilber theorem and cohomological descent for \v{C}ech nerves, we know that the map $\phi: \widehat{W}_{\bullet} \xrightarrow{\sim} V_{\bullet}$ given by sending $(u_0, \dots, u_l) \in \widehat{W}_{l}$ with $\widetilde{f_l}(u_0) = \dots = \widetilde{f_l}(u_l) = v \in V_l$ to $v$ induces a quasi-isomorphism between the $k$-truncated total complexes of differential forms on $\widehat{W}_{\bullet}$ and $V_{\bullet}$. We will use this to prove that $\widetilde{f}$ also induces a quasi-isomorphism between the $k$-truncated total complexes of differential forms on $V_{\bullet}$ and $U_{\bullet}$ by showing that it is a retract of $\phi$ resulting in the diagram
			\begin{center}
				\begin{tikzcd}
					& U_{\bullet} \arrow[shift left]{dl}{s} \arrow[twoheadrightarrow]{dr}{\widetilde{f}} & \\
					\widehat{W}_{\bullet} \arrow[shift left]{ur}{g} \arrow{rr}{\phi} & & V_{\bullet}
				\end{tikzcd}\\
			\end{center}
			with $g \circ s = \mbox{id}$. First, define $s$ levelwise as the inclusion of the diagonal copy of $U_l$ in $\widehat{W}_l$:
			\begin{equation}
				s: u\in U_l \cong W_{0l} \mapsto (s_0)^l u = (u, \dots, u) \in W_{ll} \cong \widehat{W}_l
			\end{equation}
			with $s_0$ here meaning the vertical degeneracy map $s_0: W_{ij} \rightarrow W_{(i+1)j}$. On the other hand, since $U_{\bullet} \cong \mbox{cosk}_m(X_{\bullet}/Y_{\bullet})$, the simplicial map $g$ will be fully determined by defining it levelwise up to the $m$-th level. For $0\leq l \leq m$, we will define $g$ as the projection to the last component
			\begin{equation}
				g: (u_0, \dots, u_l) \in W_{ll} \cong \widehat{W}_l \mapsto u_l\in U_l.
			\end{equation}
			This map can be written as a composition of face maps (composing until only the last face remains), thus it defines a simplicial map. Note that for $l>m$, the column $W_{\bullet l}$ is constant and hence $\widehat{W}_l \cong U_l \cong V_l$, so $g$ will be an isomorphism on these levels. Checking the compositions gives
			\begin{center}
				\begin{tikzcd}
					u \in U_l \arrow{dr}{s} \arrow{rr}{\mbox{id}} & & u\in U_l\\
					& (u, \dots, u) \in W_{ll} = \widehat{W}_l \arrow{ur}{g}&
				\end{tikzcd}\\
			\end{center}
			so $g\circ s = \mbox{id}$ holds, and
			\begin{center}
				\begin{tikzcd}
					(u_0, \dots, u_l) \in \widehat{W}_l \arrow{dr}{g} \arrow{rr}{\phi} & & \widetilde{f}(u_l) = v\in V_l\\
					& u_l \in U_l \arrow{ur}{\widetilde{f}}&
				\end{tikzcd}\\
			\end{center}
			so $\widetilde{f}\circ g = \phi$ also holds. We automatically get 
			\begin{equation}
				\phi \circ s = (\widetilde{f} \circ g) \circ s = \widetilde{f} \circ (g\circ s) = \widetilde{f} \circ \mbox{id} = \widetilde{f}.
			\end{equation} 
			This proves that $\widetilde{f}$ is a retract of $\phi$ and because $\phi^{*}$ is a quasi-isomorphism, so is $\widetilde{f}^{*}$ due to (\ref{retract}).
		\end{proof}
	\end{lemma}
	
	\subsection{Main theorem}
	
	\begin{theorem}\label{main}
		Let
		\begin{center}
			\begin{tikzcd}
				& Z_{\bullet} \arrow[swap, twoheadrightarrow]{dl}{g_{\bullet}} \arrow[twoheadrightarrow]{dr}{h_{\bullet}} & \\
				X_{\bullet} & & Y_{\bullet}
			\end{tikzcd}\\
		\end{center}
		be a Morita equivalence of Lie $n$-groupoids and $\alpha_{\bullet}$ an $m$-shifted symplectic form on $X_{\bullet}$. Then, there is also an induced $m$-shifted symplectic form $\beta_{\bullet}$ on $Y_{\bullet}$ and an $(m-1)$-shifted $2$-form $\phi_{\bullet}$ on $Z_{\bullet}$ such that the zig-zag of hypercovers becomes a symplectic Morita equivalence.
		\begin{proof}
			We now know from Lemma (\ref{descent_hypercovers}), that for a hypercover $f_{\bullet}$ between Lie $n$-groupoids, $f_{\bullet}^{*}$ is an isomorphism between the $k$-truncated total cohomology groups for any fixed $k\geq 0$, so in particular for $k=2$. Hence, we can invert $g_{\bullet}$ and $h_{\bullet}$ given above on the level of cohomology classes. Since $\alpha_{\bullet}$ is closed, we have that $[\alpha_{\bullet}] \in H_D^{2+m}(X_{\bullet}, F^2 \widehat{\Omega}^{\bullet})$. So, we can apply $g_{\bullet}^{*}$ to get 
			\begin{equation}
				g_{\bullet}^{*}[\alpha_{\bullet}] = [g_{\bullet}^{*} \alpha_{\bullet}] \in H_D^{2+m}(Z_{\bullet}, F^2 \widehat{\Omega}^{\bullet}).
			\end{equation} 
			Since $h_{\bullet}^{*}$ is an isomorphism on the level of cohomology, we can find \begin{equation}
				\beta_{\bullet} \in (h_{\bullet}^{*})^{-1} g_{\bullet}^{*}[\alpha_{\bullet}] \in H_D^{2+m}(Y_{\bullet}, F^2 \widehat{\Omega}^{\bullet})
			\end{equation} 
			with $h_{\bullet}^{*}(\beta_{\bullet}) = g_{\bullet}^{*}\alpha_{\bullet} + D\phi_{\bullet}$ for some $(m-1)$-shifted $2$-form $\phi_{\bullet}$ on $Z_{\bullet}$.
			
			Then, we get that $\beta_{\bullet}$ is a closed $m$-shifted $2$-form on $Y_{\bullet}$ by definition. Note that since $g_{\bullet}$ is a simplicial morphism, it commutes with face and degeneracy maps and thus in particular with the simplicial differential $\delta$. Additionally, pullbacks of smooth maps commute with the de Rham differential, so we have $g_{\bullet}^{*}D = Dg_{\bullet}^{*}$ and $g_{\bullet}^{*} s_i^{*} = s_i^{*} g_{\bullet}^{*}$. Similarly, $s_i^{*} g_{\bullet}^{*} \alpha_{\bullet} = g_{\bullet}^{*} s_i^{*} \alpha_{\bullet} = 0$, so $g_{\bullet}^{*} \alpha_{\bullet}$ is again normalised. Thus, we can conclude that $\beta_{\bullet}$ is normalised. 
			
			Showing that $\beta_{\bullet}$ is non-degenerate requires a bit more work. For this, we will use a similar argument to the one used to prove Lemma 2.30 in \cite{cuecazhu2023shiftedsymplectic} for strict morphisms as presented in Example (\ref{Ep_symplME}). First, note that IM-forms are invariant under gauge invariance as proven in the Appendix of \cite{cuecazhu2023shiftedsymplectic}, so
			\begin{equation}
				\lambda^{g_{\bullet}^{*}\alpha_{\bullet} + D\phi_{\bullet}} = \lambda^{g_{\bullet}^{*}\alpha_{\bullet}}.
			\end{equation} 
			
			Since non-degeneracy is defined pointwise, consider $x_0\in X_0$ and $z_0 \in Z_0$ such that $g_0(z_0)=x_0$. Such a $z_0$ exists since $g_0: Z_0 \twoheadrightarrow X_0$ is a surjective submersion. Because $g_{\bullet}$ commutes with face and degeneracy maps, we can use the formula for the IM-form given in Definition (\ref{IM_form}) and Remark (\ref{IM_properties}) to get
			\begin{equation}
				\lambda_{z_0}^{g_{\bullet}^{*} \alpha_{\bullet}}(v,w) = \lambda_{x_0}^{\alpha_{\bullet}}(T_{z_0}g(v),T_{z_0}g(w)).
			\end{equation}
			
			From Lemma (\ref{hypercoverquasi}), we get that $g_{*}: H_{\bullet}(\mathcal{T}(Z)) \xrightarrow{\sim} H_{\bullet}(\mathcal{T}(X))$ is an isomorphism, hence the pairings induced by $\lambda^{g_{\bullet}^{*} \alpha_{\bullet}}$ and $\lambda^{\alpha_{\bullet}}$ differ only by this isomorphism, so $\alpha_{\bullet}$ is non-degenerate iff $g_{\bullet}^{*} \alpha_{\bullet}$ is non-degenerate. Similarly, we can conclude that $\beta_{\bullet}$ is non-degenerate iff $h_{\bullet}^{*} \beta_{\bullet}$ is non-degenerate. Thus, in total, the non-degeneracy of $\alpha_{\bullet}$ implies that of $g_{\bullet}^{*} \alpha_{\bullet}$. This has the same IM-form as $g_{\bullet}^{*} \alpha_{\bullet} + D\phi_{\bullet} = h_{\bullet}^{*} \beta_{\bullet}$, which is hence also non-degenerate. In turn, this implies that $\beta_{\bullet}$ itself is also non-degenerate. 
			
			Lastly, to show that we get a symplectic Morita equivalence, note that
			\begin{equation}
				g_{\bullet}^{*}\alpha_{\bullet} - h_{\bullet}^{*} \beta_{\bullet} = D\phi_{\bullet}.
			\end{equation} 
			Thus, all of the properties are satisfied and the zig-zag of hypercovers becomes a symplectic Morita equivalence between $(X_{\bullet}, \alpha_{\bullet})$ and $(Y_{\bullet}, \beta_{\bullet})$.
		\end{proof}
	\end{theorem}
	
	This statement is particularly useful when dealing with different Lie $n$-groupoid models for the same $n$-stack. Considering different presentations of the same (higher) stack is often useful because depending on the context, a certain model might be more useful than a different one in one situation and \emph{vice versa} in another situation. If we now have a shifted symplectic form on one of these Lie $n$-groupoid models, we can use the theorem above to explicitly (up to some choice of gauge transformation) construct a shifted symplectic form on any one of the other models we are interested in by pulling back along the Morita equivalence between them. This will then result in compatible symplectic structures among all of the presentations for our $n$-stack in question. This can also be viewed as a symplectic structure on the corresponding $n$-stack.

	
	
	
	\appendix
	
	\section{Notes on quasi-isomorphisms}
	
	In this section, we will collect proofs for some general statements about cochain complexes and quasi-isomorphisms between them which will be needed in the proof of cohomological descent of hypercovers. All of these statements will be written in the framework of smooth manifolds as this is the framework needed for this paper. Note, however, that most of the proofs would translate to any suitable (abelian) category. Many of the statements in this section can probably be proven faster using spectral sequences, but we will not assume any knowledge on this topic and instead construct the proofs by hand.\\ 
	
	The first statement we will discuss is a generalization of the generalized Mayer-Vietoris argument from \cite{bott1982differential} to an augmented triple complex using an analogous proof. For this, consider a first octant triple complex $K^{r, p,q} \coloneqq \Omega^q(Z_{r,p})$ for a bisimplicial manifold $Z_{\bullet \bullet}$ and $r, p,q \in \Z_{\geq 0}$. Let the two horizontal differentials be the simplicial differentials $\delta_1$ in the first component and $\delta_2$ in the second component and let the vertical differential be the de Rham differential $d$ in the third component. Then, we can form the total complex $K^n \coloneqq \bigoplus_{r+p+q=n} K^{r,p,q}$ with differential $D = \delta_1 + (-1)^r \delta_2 + (-1)^{r+p} d$. Note that the following proof works analogously for the $k$-truncated triple complex $F^k \widehat{\Omega}^q(Z_{r,p})$ for fixed $k\geq 0$, as this will only lead to more terms being zero.
	
	\begin{lemma}\label{triplecx}
		Let $K^{r,p,q}$ be a triple complex with differentials as above and let $V^{p,q} \coloneqq \Omega^q(Y_p)$ be the cochain complex corresponding to differential forms on a simplicial manifold $Y_{\bullet}$ with the simplicial differential $\delta_2$ in the first component and the de Rham differential $d$ in the second component together with surjective submersions $f_p: Z_{0,p} \twoheadrightarrow Y_p$ inducing injections ${f_p}^{*}: V^{p,q} \hookrightarrow K^{0,p,q}$ for each $p,q\in \Z_{\geq 0}$ such that the rows (in $\delta_1$-direction) of the corresponding augmented triple complex are exact. Then, ${f_{\bullet}}^{*}$ applied to the total complex $V^n \coloneqq \bigoplus_{p+q=n} V^{p,q}$ with total differential $\widetilde{D} = \delta_2 + (-1)^p d$ is a quasi-isomorphism, i.e~it descends to an isomorphism between the cohomology groups of the total complex $V^n$ and those of the total complex $K^{n}$.
		\begin{proof} The proof for this statement works analogous to the argument provided in \cite{bott1982differential}, so we will start by proving that $f_{\bullet}^{*}$ is a cochain map on the total complexes and hence descends to a map between the cohomology groups. Then, we will prove surjectivity and injectivity on the level of cohomology classes. We can visualize $f_{\bullet}^{*}$ in the following way.
			\begin{center}
				\begin{tikzcd}[column sep=0.8cm]
					K^{0,0,0} \arrow{r}{D} 
					& K^{0,0,1} \oplus K^{0,1,0} \oplus K^{1,0,0} \arrow{r}{D} 
					& K^{0,0,2} \oplus K^{0,1,1} \oplus K^{0,2,0} \oplus K^{1,0,1} \oplus \dots \arrow{r}{D} 
					& \dots \\
					V^{0,0} \arrow{r}{\widetilde{D}} \arrow[hookrightarrow]{u}{f_{0}^{*}}
					& V^{0,1} \oplus V^{1,0} \textcolor{white}{\oplus 000000} \arrow[shorten <= -1.5cm]{r}{\widetilde{D}} \arrow[shift left=1.5cm, hookrightarrow]{u}{f_{0}^{*}}  \arrow[hookrightarrow]{u}{f_{1}^{*}}
					& V^{0,2} \oplus V^{1,1} \oplus V^{2,0} \textcolor{white}{\oplus 0000 \oplus 0000} \arrow[shorten <= -2.25cm]{r}{\widetilde{D}} \arrow[shift left=2.5cm, hookrightarrow]{u}{f_{0}^{*}} \arrow[shift left=1.25cm, hookrightarrow]{u}{f_{1}^{*}} \arrow[shift left=0cm, hookrightarrow]{u}{f_{2}^{*}}
					& \dots
				\end{tikzcd}
			\end{center}
			\emph{$f_{\bullet}^{*}$ is a cochain map}: The map $f_{\bullet}^{*}$ commutes with the differentials in the picture above as
			\begin{equation}
				Df_{\bullet}^{*} = (\delta_1 + \delta_2 + (-1)^p d)f_{\bullet}^{*} = \delta_1 f_{\bullet}^{*} + \delta_2 f_{\bullet}^{*} + (-1)^p df_{\bullet}^{*} = 0 + f_{\bullet}^{*}\delta_2 + f_{\bullet}^{*} d.
			\end{equation}
			This works because pullbacks commute with the de Rham differential, $f_{\bullet}^{*}$ being a simplicial map makes it commute with $\delta_2$, and exactness of the augmented rows gives $\delta_1 f_{\bullet}^{*} = 0$. Hence $f_{\bullet}^{*}$ is a cochain map and descends to a well-defined map on the level of cohomology.\\
			
			\emph{$f_{\bullet}^{*}$ is surjective on cohomology}: Let $\phi \in K^{n}$ with $D\phi = 0$. Then, $\delta_1 \phi_k = 0$ for all of the non-zero components $\phi_k \in K^{r, l, k-l}$ of $\phi$ with $r$ being the maximal index in $\delta_1$-direction such that $\phi$ contains non-zero components on that level. By exactness of the augmented rows in $\delta_1$-direction, we have $\phi_k = \delta_1 \psi_k$ for some $\psi_k \in \Omega^{k-l}(X_{(r-1), l})$. Therefore, $\widetilde{\phi} \coloneqq \phi - D\psi_k$ is in the same cohomology class as $\phi$ and has one less non-zero component at index $r$ in $\delta_1$-direction (see Fig.\ref{fig2}). Iterating this procedure, we arrive at $\widetilde{\phi} = \phi + D\psi$ with the only non-zero components of $\widetilde{\phi}$ lying in the $r=0$ layer of $K^n$. Therefore, by exactness of the rows in $\delta_1$-direction, all components of $\widetilde{\phi}$ lie in the image of ${f_{\bullet}}^{*}$ and hence $\widetilde{\phi} \in \hbox{im}(f_{\bullet}^{*})$.
			
			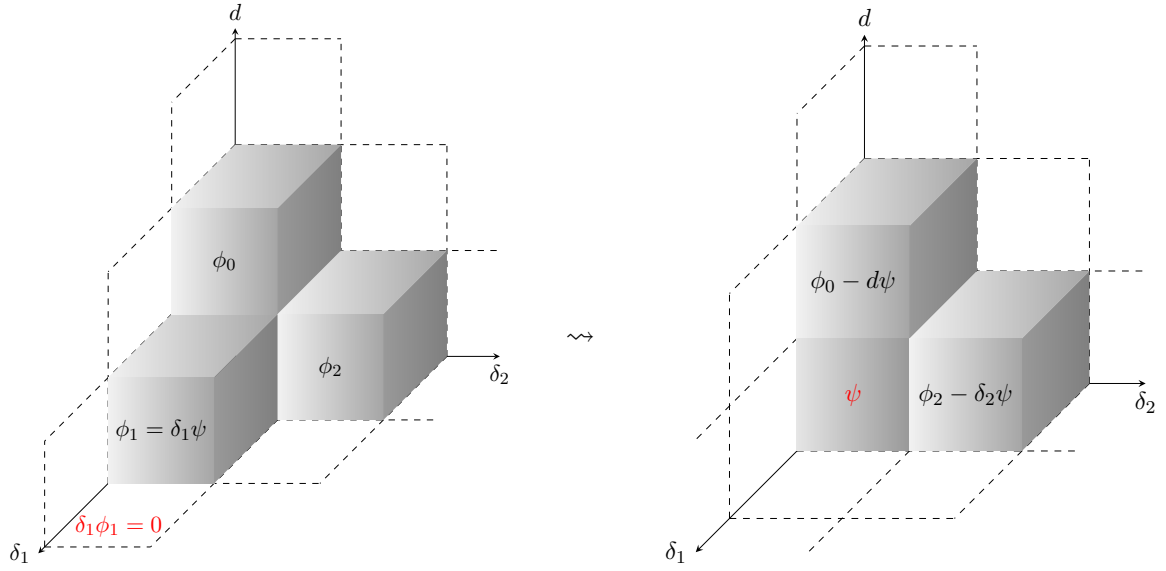
\begin{figure}[htbp]
				\centering
				
				\begin{subfigure}[b]{0.44\textwidth}
					\centering
					\scalebox{.8}{
						\begin{tikzpicture}[x=1.75cm,y=1.75cm,z=-1.05cm,>=stealth]
							
							\draw[->] (xyz cs:x=0) -- (xyz cs:x=2.5) node[below] {$\delta_2$};
							\draw[->] (xyz cs:y=0) -- (xyz cs:y=3.1) node[above] {$d$};
							\draw[->] (xyz cs:z=0) -- (xyz cs:z=3.1) node[left] {$\delta_1$};
							
							\draw[dashed] (1,0) -- (1,3);
							\draw[dashed] (2,0) -- (2,2);
							
							\draw[dashed] (0,1) -- (2.5,1);
							\draw[dashed] (0,2) -- (2,2);
							\draw[dashed] (0,3) -- (1,3);
							
							\draw[dashed] (xyz cs: z=1) -- (xyz cs: x=2.5,z=1);
							\draw[dashed] (xyz cs: x=1) -- (xyz cs: x=1,z=3);
							
							\draw[dashed] (xyz cs: z=2) -- (xyz cs: x=2,z=2);
							\draw[dashed] (xyz cs: x=2) -- (xyz cs: x=2,z=2);
							
							\draw[dashed] (xyz cs: z=3) -- (xyz cs: x=1,z=3);
							
							\draw[dashed] (xyz cs: z=1) -- (xyz cs: y=3,z=1);
							\draw[dashed] (xyz cs: y=1) -- (xyz cs: y=1,z=3);
							
							\draw[dashed] (xyz cs: z=2) -- (xyz cs: y=2,z=2);
							\draw[dashed] (xyz cs: y=2) -- (xyz cs: y=2,z=2);
							
							\draw[dashed] (xyz cs: z=3) -- (xyz cs: y=1,z=3);
							\draw[dashed] (xyz cs: y=3) -- (xyz cs: y=3,z=1);
							
							\begin{scope}[every node/.append style={
									xslant=1},xslant=1]
								\shade[right color=gray!70,left color=gray!10] (xyz cs: x=-1.4,y=2,z=1) rectangle (xyz cs: x=-1,y=2,z=0);
								\shade[right color=gray!70,left color=gray!10] (xyz cs: x=0.2,y=1,z=2) rectangle (xyz cs: x=0.6,y=1,z=1);
								\shade[right color=gray!70,left color=gray!10] (xyz cs: x=0.6,y=1,z=1) rectangle (xyz cs: x=1,y=1,z=0);
							\end{scope}
							
							\begin{scope}
								\shade[left color=gray!10, right color=gray!70] (xyz cs: x=0,y=2,z=1) rectangle (xyz cs: x=1,y=1,z=1);
								\shade[left color=gray!10, right color=gray!70] (xyz cs: x=0,y=1,z=2) rectangle (xyz cs: x=1,y=0,z=2);
								\shade[left color=gray!10, right color=gray!70] (xyz cs: x=1,y=1,z=1) rectangle (xyz cs: x=2,y=0,z=1);
							\end{scope}
							
							\begin{scope}[every node/.append style={
									yslant=1},yslant=1]
								\shade[left color=gray!80,right color=black!50] (xyz cs: x=1,y=0.6,z=1) rectangle (xyz cs: x=1,y=1,z=0);
								\shade[left color=gray!80,right color=black!50] (xyz cs: x=1,y=0.2,z=2) rectangle (xyz cs: x=1,y=0.6,z=1);
								\shade[left color=gray!80,right color=black!50] (xyz cs: x=2,y=-1.4,z=1) rectangle (xyz cs: x=2,y=-1,z=0);
							\end{scope}
							
							\node at (xyz cs: x=0.5,y=1.5,z=1) {$\phi_0$};
							\node at (xyz cs: x=0.5,y=0.5,z=2) {$\phi_1 = \delta_1 \psi$};
							\node at (xyz cs: x=1.5,y=0.5,z=1) {$\phi_2$};
							
							\node at (xyz cs: x=1,y=0.5,z=3.5) {\textcolor{red}{$\delta_1 \phi_1 = 0$}};
						\end{tikzpicture}
					}
				\end{subfigure}
				\begin{subfigure}[b]{.1\textwidth}
					\centering
					$\rightsquigarrow$\\
					\vspace{3cm}
				\end{subfigure}
				\begin{subfigure}[b]{0.44\textwidth}
					\centering
					\scalebox{.85}{
						\begin{tikzpicture}[x=1.75cm,y=1.75cm,z=-1.05cm,>=stealth]
							
							\draw[->] (xyz cs:x=0) -- (xyz cs:x=2.5) node[below] {$\delta_2$};
							\draw[->] (xyz cs:y=0) -- (xyz cs:y=3.1) node[above] {$d$};
							\draw[->] (xyz cs:z=0) -- (xyz cs:z=2.5) node[left] {$\delta_1$};
							
							\draw[dashed] (1,0) -- (1,3);
							\draw[dashed] (2,0) -- (2,2);
							
							\draw[dashed] (0,1) -- (2.5,1);
							\draw[dashed] (0,2) -- (2,2);
							\draw[dashed] (0,3) -- (1,3);
							
							\draw[dashed] (xyz cs: z=1) -- (xyz cs: x=2.5,z=1);
							\draw[dashed] (xyz cs: x=1) -- (xyz cs: x=1,z=2.5);
							
							\draw[dashed] (xyz cs: z=2) -- (xyz cs: x=2,z=2);
							\draw[dashed] (xyz cs: x=2) -- (xyz cs: x=2,z=2);
							
							
							\draw[dashed] (xyz cs: z=1) -- (xyz cs: y=3,z=1);
							\draw[dashed] (xyz cs: y=1) -- (xyz cs: y=1,z=2.5);
							
							\draw[dashed] (xyz cs: z=2) -- (xyz cs: y=2,z=2);
							\draw[dashed] (xyz cs: y=2) -- (xyz cs: y=2,z=2);
							
							\draw[dashed] (xyz cs: y=3) -- (xyz cs: y=3,z=1);
							
							\begin{scope}[every node/.append style={
									xslant=1},xslant=1]
								\shade[right color=gray!70,left color=gray!10] (xyz cs: x=-1.4,y=2,z=1) rectangle (xyz cs: x=-1,y=2,z=0);
								\shade[right color=gray!70,left color=gray!10] (xyz cs: x=0.6,y=1,z=1) rectangle (xyz cs: x=1,y=1,z=0);
							\end{scope}
							
							\begin{scope}
								\shade[left color=gray!10, right color=gray!70] (xyz cs: x=0,y=2,z=1) rectangle (xyz cs: x=1,y=1,z=1);
								\shade[left color=gray!20, right color=gray!80] (xyz cs: x=0,y=0,z=1) rectangle (xyz cs: x=1,y=1,z=1);
								\shade[left color=gray!10, right color=gray!70] (xyz cs: x=1,y=1,z=1) rectangle (xyz cs: x=2,y=0,z=1);
							\end{scope}
							
							\begin{scope}[every node/.append style={
									yslant=1},yslant=1]
								\shade[left color=gray!80,right color=black!50] (xyz cs: x=1,y=0.6,z=1) rectangle (xyz cs: x=1,y=1,z=0);
								\shade[left color=gray!80,right color=black!50] (xyz cs: x=2,y=-1.4,z=1) rectangle (xyz cs: x=2,y=-1,z=0);
							\end{scope}
							
							\node at (xyz cs: x=0.5,y=1.5,z=1) {$\phi_0- d\psi$};
							\node at (xyz cs: x=1.5,y=0.5,z=1) {$\phi_2 - \delta_2 \psi$};
							
							\node at (xyz cs: x=0.5,y=0.5,z=1) {\textcolor{red}{$\psi$}};
							
						\end{tikzpicture}
					}
				\end{subfigure}
				
				\caption{Surjectivity on cohomology of $f^{*}_{\bullet}$.}
				\label{fig2}
			\end{figure}
			
			\emph{$f_{\bullet}^{*}$ is injective on cohomology}: Let $\omega = \omega_0 + \dots + \omega_n \in V^n$ with $\omega_k \in \Omega^{n-k}(Y_k)$. If $f_{\bullet}^{*}(\omega) \eqqcolon \psi = \psi_0 + \dots + \psi_n \in K^n$ with $\psi_k \in \Omega^{n-k}(Z_{0,k})$ is exact, i.e.~if there is some $\phi \in K^{n-1}$ with $f_{\bullet}^{*}(\omega) = D\phi$, then $\delta_1 \phi_i = 0$ for all $i$ (in particular for the last non-zero component). Thus, we can again construct some $\widetilde{\phi} \in K^{n-1}$ in the same cohomology class which has its only non-zero components in the $(r=0)$-layer (see Fig.\ref{fig3}). This means $\widetilde{\phi} = \widetilde{\phi}_0 + \dots + \widetilde{\phi}_{n-1}$ with $\widetilde{\phi}_k \in \Omega^{n-k-1}(Z_{0,k})$. Then, because $\delta_1 \widetilde{\phi}_k = 0$ for all components in $\widetilde{\phi}$, we have $f_{\bullet}^{*}(\omega) = D\phi = D\widetilde{\phi} = (\delta_2 + (-1)^p d)\widetilde{\phi}$. Due to the exactness of the augmented rows, we can again conclude $\widetilde{\phi} =f_{\bullet}^{*}(\omega')$ for some $\omega' \in V^{n-1}$ with $\widetilde{D}\omega' = \omega$. Hence, $\omega$ is also exact.
			
			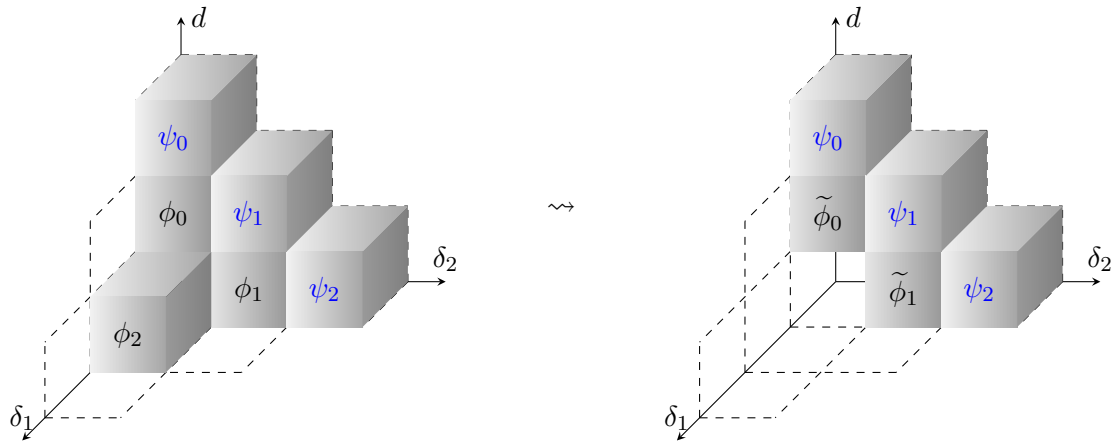
\begin{figure}[htbp]
				\centering
				
				\begin{subfigure}[b]{0.44\textwidth}
					\centering
					
					\begin{tikzpicture}[x=1cm,y=1cm,z=-0.6cm,>=stealth]
						
						\draw[->] (xyz cs:x=0) -- (xyz cs:x=3.5) node[above] {$\delta_2$};
						\draw[->] (xyz cs:y=0) -- (xyz cs:y=3.5) node[right] {$d$};
						\draw[->] (xyz cs:z=0) -- (xyz cs:z=3.5) node[above] {$\delta_1$};
						
						\draw[dashed] (1,0) -- (1,3);
						\draw[dashed] (2,0) -- (2,2);
						\draw[dashed] (3,0) -- (3,1);
						
						\draw[dashed] (0,1) -- (3,1);
						\draw[dashed] (0,2) -- (2,2);
						\draw[dashed] (0,3) -- (1,3);
						
						\draw[dashed] (xyz cs: z=1) -- (xyz cs: x=3,z=1);
						\draw[dashed] (xyz cs: x=1) -- (xyz cs: x=1,z=3);
						
						\draw[dashed] (xyz cs: z=2) -- (xyz cs: x=2,z=2);
						\draw[dashed] (xyz cs: x=2) -- (xyz cs: x=2,z=2);
						
						\draw[dashed] (xyz cs: z=3) -- (xyz cs: x=1,z=3);
						\draw[dashed] (xyz cs: x=3) -- (xyz cs: x=3,z=1);
						
						\draw[dashed] (xyz cs: z=1) -- (xyz cs: y=3,z=1);
						\draw[dashed] (xyz cs: y=1) -- (xyz cs: y=1,z=3);
						
						\draw[dashed] (xyz cs: z=2) -- (xyz cs: y=2,z=2);
						\draw[dashed] (xyz cs: y=2) -- (xyz cs: y=2,z=2);
						
						\draw[dashed] (xyz cs: z=3) -- (xyz cs: y=1,z=3);
						\draw[dashed] (xyz cs: y=3) -- (xyz cs: y=3,z=1);
						
						
						\begin{scope}[every node/.append style={
								xslant=1},xslant=1]
							\shade[right color=gray!70,left color=gray!10] (xyz cs: x=0.2,y=1,z=2) rectangle (xyz cs: x=0.6,y=1,z=1);
						\end{scope}
						
						\begin{scope}
							\shade[left color=gray!20, right color=gray!80] (xyz cs: x=0,y=2,z=1) rectangle (xyz cs: x=1,y=1,z=1);
							\shade[left color=gray!10, right color=gray!70] (xyz cs: x=0,y=1,z=2) rectangle (xyz cs: x=1,y=0,z=2);
							\shade[left color=gray!20, right color=gray!80] (xyz cs: x=1,y=1,z=1) rectangle (xyz cs: x=2,y=0,z=1);
						\end{scope}
						
						\begin{scope}[every node/.append style={
								yslant=1},yslant=1]
							\shade[left color=gray!80,right color=black!50] (xyz cs: x=1,y=0.2,z=2) rectangle (xyz cs: x=1,y=0.6,z=1);
						\end{scope}

						
						\begin{scope}[every node/.append style={
								xslant=1},xslant=1]
							\shade[right color=gray!70,left color=gray!10] (xyz cs: x=-2.4,y=3,z=1) rectangle (xyz cs: x=-2,y=3,z=0);
							\shade[right color=gray!70,left color=gray!10] (xyz cs: x=1.6,y=1,z=1) rectangle (xyz cs: x=2,y=1,z=0);
							\shade[right color=gray!70,left color=gray!10] (xyz cs: x=-0.4,y=2,z=1) rectangle (xyz cs: x=0,y=2,z=0);
						\end{scope}
						
						\begin{scope}
							\shade[left color=gray!10, right color=gray!70] (xyz cs: x=0,y=2,z=1) rectangle (xyz cs: x=1,y=3,z=1);
							\shade[left color=gray!10, right color=gray!70] (xyz cs: x=2,y=0,z=1) rectangle (xyz cs: x=3,y=1,z=1);
							\shade[left color=gray!10, right color=gray!70] (xyz cs: x=1,y=2,z=1) rectangle (xyz cs: x=2,y=1,z=1);
						\end{scope}
						
						\begin{scope}[every node/.append style={
								yslant=1},yslant=1]
							\shade[left color=gray!80,right color=black!50] (xyz cs: x=1,y=1.6,z=1) rectangle (xyz cs: x=1,y=2,z=0);
							\shade[left color=gray!80,right color=black!50] (xyz cs: x=3,y=-2.4,z=1) rectangle (xyz cs: x=3,y=-2,z=0);
							\shade[left color=gray!80,right color=black!50] (xyz cs: x=2,y=-0.4,z=1) rectangle (xyz cs: x=2,y=0,z=0);
						\end{scope}
						
						\node at (xyz cs: x=0.5,y=2.5,z=1) {\textcolor{blue}{$\psi_0$}};
						\node at (xyz cs: x=1.5,y=1.5,z=1) {\textcolor{blue}{$\psi_1$}};
						\node at (xyz cs: x=2.5,y=0.5,z=1) {\textcolor{blue}{$\psi_2$}};
						
						\node at (xyz cs: x=0.5,y=1.5,z=1) {$\phi_0$};
						\node at (xyz cs: x=1.5,y=0.5,z=1) {$\phi_1$};
						\node at (xyz cs: x=0.5,y=0.5,z=2) {$\phi_2$};
						
					\end{tikzpicture}
				\end{subfigure}
				\begin{subfigure}[b]{.1\textwidth}
					\centering
					$\rightsquigarrow$\\
					\vspace{3cm}
				\end{subfigure}
				\begin{subfigure}[b]{0.44\textwidth}
					\centering
					
					\begin{tikzpicture}[x=1cm,y=1cm,z=-0.6cm,>=stealth]
						
						\draw[->] (xyz cs:x=0) -- (xyz cs:x=3.5) node[above] {$\delta_2$};
						\draw[->] (xyz cs:y=0) -- (xyz cs:y=3.5) node[right] {$d$};
						\draw[->] (xyz cs:z=0) -- (xyz cs:z=3.5) node[above] {$\delta_1$};
						
						\draw[dashed] (1,0) -- (1,3);
						\draw[dashed] (2,0) -- (2,2);
						\draw[dashed] (3,0) -- (3,1);
						
						\draw[dashed] (0,1) -- (3,1);
						\draw[dashed] (0,2) -- (2,2);
						\draw[dashed] (0,3) -- (1,3);
						
						\draw[dashed] (xyz cs: z=1) -- (xyz cs: x=3,z=1);
						\draw[dashed] (xyz cs: x=1) -- (xyz cs: x=1,z=3);
						
						\draw[dashed] (xyz cs: z=2) -- (xyz cs: x=2,z=2);
						\draw[dashed] (xyz cs: x=2) -- (xyz cs: x=2,z=2);
						
						\draw[dashed] (xyz cs: z=3) -- (xyz cs: x=1,z=3);
						\draw[dashed] (xyz cs: x=3) -- (xyz cs: x=3,z=1);
						
						\draw[dashed] (xyz cs: z=1) -- (xyz cs: y=3,z=1);
						\draw[dashed] (xyz cs: y=1) -- (xyz cs: y=1,z=3);
						
						\draw[dashed] (xyz cs: z=2) -- (xyz cs: y=2,z=2);
						\draw[dashed] (xyz cs: y=2) -- (xyz cs: y=2,z=2);
						
						\draw[dashed] (xyz cs: z=3) -- (xyz cs: y=1,z=3);
						\draw[dashed] (xyz cs: y=3) -- (xyz cs: y=3,z=1);
						
						
						\begin{scope}[every node/.append style={
								xslant=1},xslant=1]
						\end{scope}
						
						\begin{scope}
							\shade[left color=gray!20, right color=gray!80] (xyz cs: x=0,y=2,z=1) rectangle (xyz cs: x=1,y=1,z=1);
							\shade[left color=gray!20, right color=gray!80] (xyz cs: x=1,y=1,z=1) rectangle (xyz cs: x=2,y=0,z=1);
						\end{scope}
						
						\begin{scope}[every node/.append style={
								yslant=1},yslant=1]
						\end{scope}

						
						\begin{scope}[every node/.append style={
								xslant=1},xslant=1]
							\shade[right color=gray!70,left color=gray!10] (xyz cs: x=-2.4,y=3,z=1) rectangle (xyz cs: x=-2,y=3,z=0);
							\shade[right color=gray!70,left color=gray!10] (xyz cs: x=1.6,y=1,z=1) rectangle (xyz cs: x=2,y=1,z=0);
							\shade[right color=gray!70,left color=gray!10] (xyz cs: x=-0.4,y=2,z=1) rectangle (xyz cs: x=0,y=2,z=0);
						\end{scope}
						
						\begin{scope}
							\shade[left color=gray!10, right color=gray!70] (xyz cs: x=0,y=2,z=1) rectangle (xyz cs: x=1,y=3,z=1);
							\shade[left color=gray!10, right color=gray!70] (xyz cs: x=2,y=0,z=1) rectangle (xyz cs: x=3,y=1,z=1);
							\shade[left color=gray!10, right color=gray!70] (xyz cs: x=1,y=2,z=1) rectangle (xyz cs: x=2,y=1,z=1);
						\end{scope}
						
						\begin{scope}[every node/.append style={
								yslant=1},yslant=1]
							\shade[left color=gray!80,right color=black!50] (xyz cs: x=1,y=1.6,z=1) rectangle (xyz cs: x=1,y=2,z=0);
							\shade[left color=gray!80,right color=black!50] (xyz cs: x=3,y=-2.4,z=1) rectangle (xyz cs: x=3,y=-2,z=0);
							\shade[left color=gray!80,right color=black!50] (xyz cs: x=2,y=-0.4,z=1) rectangle (xyz cs: x=2,y=0,z=0);
						\end{scope}
						
						\node at (xyz cs: x=0.5,y=2.5,z=1) {\textcolor{blue}{$\psi_0$}};
						\node at (xyz cs: x=1.5,y=1.5,z=1) {\textcolor{blue}{$\psi_1$}};
						\node at (xyz cs: x=2.5,y=0.5,z=1) {\textcolor{blue}{$\psi_2$}};
						
						\node at (xyz cs: x=0.5,y=1.5,z=1) {$\widetilde{\phi}_0$};
						\node at (xyz cs: x=1.5,y=0.5,z=1) {$\widetilde{\phi}_1$};
						
					\end{tikzpicture}
					
				\end{subfigure}
				
				\caption{Injectivity on cohomology of $f^{*}_{\bullet}$.}
				\label{fig3}
			\end{figure}
		\end{proof}
	\end{lemma}
	
	Besides these direct proofs, we require two other approaches to proving that a given map between cochain complexes is a quasi-isomorphism or that a given cochain complex is contractible in some way. The first one uses the existence of a homotopy operator and the second one is a statement about retracts of quasi-isomorphisms.
	
	\begin{lemma}\label{homotopyop}
		Let $(C^{\bullet}, d)$ be a cochain complex with differential $d: C^{k} \rightarrow C^{k+1}$ and let $K: C^{k+1} \rightarrow C^{k}$ be a \emph{homotopy operator} on $C^{\bullet}$, i.e.~$K\circ d + d\circ K = \mbox{id}$. Then, the complex $C^{\bullet}$ is \emph{contractible}, i.e.~$H_{d}^{k}(C) = 0$ for $k>0$.
		\begin{proof}
			This statement is proven in a more general sense in chapter 4 of \cite{bott1982differential} and we will adapt the proof to the specific situation needed here. Let $\omega \in C^{k}$ be closed, i.e.~$d\omega=0$. Then,  
			\begin{equation}
				\omega = \mbox{id}(\omega) = (K\circ d + d\circ K)(\omega) = (d \circ K)(\omega) + 0 = d(K\omega)
			\end{equation}
			is exact, so the identity map is the zero map in cohomology. Hence, all cohomology groups $H_{d}^{k}(C)$ with $k>0$ have to vanish.
		\end{proof}
	\end{lemma}
	
	\begin{remark}\label{homotopyrmk}
		If $(C^{\bullet}, d)$ is the total complex of a double complex which can be augmented by an injection as before, suppose that the horizontal differential admits a homotopy operator. Then, by this lemma, we have that the rows of the augmented double complex are exact and hence the total cohomology of the double complex is isomorphic to the cohomology of the augmentation.
	\end{remark}
	
	
	\begin{lemma}\label{retract}
		Let $A,B,C$ be (simplicial) manifolds. Let $\phi: C\rightarrow B$ be a morphism such that $\phi^{*}: \Omega^{\bullet}(B) \rightarrow \Omega^{\bullet}(C)$ is a quasi-isomorphism and let $f: A\rightarrow B$ be a retract of $\phi$, i.e.~there exist morphisms $g: C \rightarrow A$ and $s: A \rightarrow C$ such that $g\circ s = \mbox{id}$ and the following diagram commutes.
		\begin{center}
			\begin{tikzcd}
				& A \arrow[shift left]{dl}{s} \arrow[twoheadrightarrow]{dr}{f} & \\
				C \arrow[shift left]{ur}{g} \arrow{rr}{\phi} & & B
			\end{tikzcd}\\
		\end{center}
		Then, $f^{*}: \Omega^{\bullet}(B) \rightarrow \Omega^{\bullet}(A)$ is also a quasi-isomorphism.
		\begin{proof}
			Consider the following diagram on the level of cohomology
			\begin{center}
				\begin{tikzcd}
					H^{k}(A, \Omega^{\bullet}) \arrow{r}{g^{*}} \arrow[bend left=25]{rr}{\mbox{id}}
					& H^{k}(C, \Omega^{\bullet}) \arrow{r}{s^{*}} 
					& H^{k}(A, \Omega^{\bullet}) \\
					H^{k}(B, \Omega^{\bullet})  \arrow{r}{\mbox{id}} \arrow[swap]{r}{\cong}  \arrow{u}{f^{*}}
					& H^{k}(B, \Omega^{\bullet})  \arrow{r}{\mbox{id}} \arrow[swap]{r}{\cong}  \arrow{u}{\phi^{*}} \arrow[swap]{u}{\cong} 
					& H^{k}(B, \Omega^{\bullet}) \arrow{u}{f^{*}} \\
				\end{tikzcd}\\
			\end{center}
			From this diagram, we can see that $s^{*} \circ g^{*} = (g\circ s)^{*} = \mbox{id}$, hence $s^{*}$ must be surjective and $g^{*}$ must be injective. With this, we can conclude
			\begin{enumerate}
				\item $f^{*}$ is injective since $g^{*} \circ f^{*} = \phi^{*}$ is injective,
				\item $f^{*}$ is surjective since $f^{*} = s^{*} \circ \phi^{*}$ is surjective.
			\end{enumerate}
			This proves that $f^{*}: H^{k}(B, \Omega^{\bullet}) \xrightarrow{\sim} H^{k}(A, \Omega^{\bullet})$ is an isomorphism.
		\end{proof}
	\end{lemma}
	
	\begin{remark}
		Note that the de Rham functor $\Omega^{\bullet}$ is simply a placeholder here, as the same argument works for any cochain complexes on $A, B$ and $C$, as long as the morphisms between the objects induce cochain maps via their pullbacks. In particular, the same argument works for the $k$-truncated (and normalised) de Rham functor $F^k \widehat{\Omega}^{\bullet}$ for any $k\geq 0$.
	\end{remark}

	
	\section*{Acknowledgement}
	
	I would like to thank all participants of the BGW 2024 PhD Retreat for many inspiring discussions and helpful suggestions during my work on this project. In particular, I would like to thank my advisor, Professor Chenchang Zhu, for her guidance and constructive feedback. Additionally, I would like to thank Miquel Cueca and Florian Dorsch for their useful comments on an earlier version of this paper.
	
	\printbibliography

@article{dold1961homologie,
	author = {Dold, Albrecht and Puppe, Dieter},
	journal = {Annales de l'institut {F}ourier},
	pages = {201--312},
	title = {Homologie nicht-additiver {F}unktoren. {A}nwendungen},
	volume = {11},
	year = {1961}}

@article{kosmann2016multiplicativity,
	author = {Kosmann-Schwarzbach, Yvette},
	journal = {Banach Center Publications},
	pages = {131--166},
	title = {Multiplicativity, from {L}ie groups to generalized geometry},
	volume = {110},
	year = {2016}}

@article{grothendieck1961techniques,
	author = {Grothendieck, Alexander},
	journal = {S{\'e}minaire Bourbaki},
	number = {212},
	pages = {99--118},
	title = {Techniques de construction et th{\'e}or{\`e}mes d'existence en g{\'e}om{\'e}trie alg{\'e}brique {III}: pr{\'e}sch{\'e}mas quotients},
	volume = {6},
	year = {1961}}

@book{bates1997lectures,
	author = {Bates, Sean and Weinstein, Alan},
	publisher = {American Mathematical Soc.},
	title = {Lectures on the Geometry of Quantization},
	volume = {8},
	year = {1997}}

@article{alfonsi2023higher,
	author = {Alfonsi, Luigi},
	journal = {Preprint arXiv:2312.07308},
	title = {Higher geometry in physics},
	year = {2023}}

@article{noether1918invariante,
	author = {Noether, Emmy},
	journal = {Nachrichten von der K{\"o}niglichen Gesellschaft der Wissenschaften zu G{\"o}ttingen, Mathematisch-physikalische Klasse},
	pages = {235--257},
	title = {Invariante {V}ariationsprobleme},
	volume = {2},
	year = {1918}}

@inproceedings{duskin1979higher,
	author = {Duskin, John W.},
	booktitle = {Applications of {S}heaves ({P}roc. {R}es. {S}ympos. {A}ppl. {S}heaf {T}heory to {L}ogic, {A}lgebra and {A}nal., {U}niv. {D}urham, {D}urham, 1977)},
	editor = {Fourman, M. P. and Mulvey, C. J. and Scott, D. S.},
	pages = {255--279},
	publisher = {Springer},
	series = {Lecture {N}otes in {M}ath.},
	title = {Higher dimensional torsors and the cohomology of topoi: the abelian theory},
	volume = {753},
	year = {1979}}

@article{noohi2005foundations,
	author = {Noohi, Behrang},
	journal = {Preprint arXiv:math/0503247},
	title = {Foundations of topological stacks {I}},
	year = {2005}}

@article{ASENS_1963_3_80_4_349_0,
	author = {Ehresmann, Charles},
	journal = {Annales scientifiques de l'{\'E}cole Normale Sup{\'e}rieure},
	number = {4},
	pages = {349--426},
	publisher = {Elsevier},
	title = {Cat{\'e}gories structur{\'e}es},
	volume = {3e s{\'e}rie, 80},
	year = {1963}}

@article{blohmann2008stacky,
	author = {Blohmann, Christian},
	journal = {International Mathematics Research Notices},
	pages = {Article rnn082},
	publisher = {Oxford University Press},
	title = {Stacky {L}ie groups},
	volume = {2008},
	year = {2008}}

@article{wolfson2016descent,
	author = {Wolfson, Jesse},
	journal = {Advances in Mathematics},
	pages = {527--575},
	publisher = {Elsevier},
	title = {Descent for n-bundles},
	volume = {288},
	year = {2016}}

@article{del2024cohomology,
	author = {del Hoyo, Matias and Ortiz, Cristian and Studzinski, Fernando},
	journal = {Preprint arXiv:2410.02570},
	title = {On the cohomology of differentiable stacks},
	year = {2024}}

@article{weinstein1987symplectic,
	author = {Weinstein, Alan},
	journal = {Bulletin of the American Mathematical Society},
	number = {1},
	pages = {101--104},
	title = {Symplectic groupoids and {P}oisson manifolds},
	volume = {16},
	year = {1987}}

@article{li2023differentiating,
	author = {Li, Du and Ryvkin, Leonid and Wessel, Arne and Zhu, Chenchang},
	journal = {Preprint arXiv:2309.00901},
	title = {Differentiating ${L}_{\infty}$ groupoids. Part {I}},
	year = {2023}}

@article{cattaneo2004integration,
	author = {Cattaneo, Alberto S. and Xu, Ping},
	journal = {Journal of Geometry and Physics},
	number = {2},
	pages = {187--196},
	publisher = {Elsevier},
	title = {Integration of twisted Poisson structures},
	volume = {49},
	year = {2004}}

@article{BCWZ2004dirac,
	author = {Henrique Bursztyn and Marius Crainic and Alan Weinstein and Chenchang Zhu},
	journal = {Duke Mathematical Journal},
	number = {3},
	pages = {549--607},
	publisher = {Duke University Press},
	title = {{Integration of twisted Dirac brackets}},
	volume = {123},
	year = {2004}}

@article{bursztyn2012multiplicative,
	author = {Bursztyn, Henrique and Cabrera, Alejandro},
	journal = {Mathematische Annalen},
	pages = {663--705},
	publisher = {Springer},
	title = {Multiplicative forms at the infinitesimal level},
	volume = {353},
	year = {2012}}

@article{bursztyn2009linear,
	author = {Bursztyn, Henrique and Cabrera, Alejandro and Ortiz, Cristi{\'a}n},
	journal = {Letters in Mathematical Physics},
	pages = {59--83},
	publisher = {Springer},
	title = {Linear and multiplicative 2-forms},
	volume = {90},
	year = {2009}}

@article{eilenberg1954groups,
	author = {Eilenberg, Samuel and MacLane, Saunders},
	journal = {Annals of Mathematics},
	number = {1},
	pages = {49--139},
	publisher = {JSTOR},
	title = {On the Groups {H}($\Pi$, n), {II}: Methods of Computation},
	volume = {60},
	year = {1954}}

@article{duskin2002simplicial,
	author = {Duskin, John W.},
	journal = {Theory and Applications of Categories},
	number = {10},
	pages = {198--308},
	title = {Simplicial matrices and the nerves of weak n-categories {I}: Nerves of bicategories},
	volume = {9},
	year = {2002}}

@article{dugger_hollander_isaksen_2004,
	author = {Dugger, Daniel and Hollander, Sharon and Isaksen, Daniel C.},
	doi = {10.1017/S0305004103007175},
	journal = {Mathematical Proceedings of the Cambridge Philosophical Society},
	number = {1},
	pages = {9--51},
	title = {Hypercovers and simplicial presheaves},
	volume = {136},
	year = {2004}}

@article{xu2004momentum,
	author = {Xu, Ping},
	journal = {Journal of Differential Geometry},
	number = {2},
	pages = {289--333},
	publisher = {Lehigh University},
	title = {Momentum maps and {M}orita equivalence},
	volume = {67},
	year = {2004}}

@article{behrend2004cohomology,
	author = {Behrend, Kai},
	journal = {Intersection theory and moduli, ICTP Lect. Notes},
	pages = {249--294},
	title = {Cohomology of stacks},
	volume = {19},
	year = {2004}}

@article{del2020morita,
	author = {del Hoyo, Matias and Ortiz, Cristian},
	journal = {International Mathematics Research Notices},
	number = {14},
	pages = {4395--4432},
	publisher = {Oxford University Press},
	title = {Morita equivalences of vector bundles},
	volume = {2020},
	year = {2020}}

@book{bott1982differential,
	author = {Bott, Raoul and Tu, Loring W.},
	publisher = {Springer},
	title = {Differential forms in algebraic topology},
	volume = {82},
	year = {1982}}

@article{PTVV2013shifted,
	author = {Pantev, Tony and To{\"e}n, Bertrand and Vaqui{\'e}, Michel and Vezzosi, Gabriele},
	journal = {Publications math{\'e}matiques de l'IH{\'E}S},
	pages = {271--328},
	title = {Shifted symplectic structures},
	volume = {117},
	year = {2013}}

@misc{getzlerslides,
	author = {Ezra Getzler},
	howpublished = {\url{https://sites.northwestern.edu/getzler/}},
	title = {Differential forms on stacks (slides).},
	year = {(2014)}}

@article{henriques2008integrating,
	author = {Henriques, Andr{\'e}},
	journal = {Compositio Mathematica},
	number = {4},
	pages = {1017--1045},
	publisher = {London Mathematical Society},
	title = {Integrating ${L}_{\infty}$-algebras},
	volume = {144},
	year = {2008}}

@article{behrend2011diffstacksgerbes,
	author = {Behrend, Kai and Xu, Ping},
	journal = {Journal of Symplectic Geometry},
	number = {3},
	pages = {285--341},
	publisher = {International Press of Boston},
	title = {Differentiable Stacks and Gerbes},
	volume = {9},
	year = {2011}}

@article{zhu2009stacky,
	author = {Zhu, Chenchang},
	doi = {10.1093/imrn/rnp080},
	journal = {International Mathematics Research Notices},
	number = {21},
	pages = {4087-4141},
	title = {n-Groupoids and Stacky Groupoids},
	volume = {2009},
	year = {2009}}

@article{cuecazhu2023shiftedsymplectic,
	author = {Miquel Cueca and Chenchang Zhu},
	journal = {Advances in Mathematics},
	pages = {Article 108829},
	title = {Shifted symplectic higher {L}ie groupoids and classifying spaces},
	volume = {413},
	year = {2023}}
	
\end{document}